\documentclass[12pt]{article}
\usepackage{amsfonts}

\usepackage{latexsym}
\usepackage{amssymb}
\usepackage{epsfig}
\usepackage{amsmath}
\usepackage{amsthm}
\usepackage[mathscr]{eucal}
\usepackage{subfigure}
\usepackage{graphicx}

\newtheorem{Theorem}{Theorem}

\newtheorem{Lemma}{Lemma}[section]

\renewcommand{\qed}{\hfill{\ \ \rule{2mm}{2mm}} \vspace{0.2in}}

\newcommand{\ind}{1\hspace{-2.3mm}{1}}

\begin{document}

\title{Phase transition in inhomogenous Erd\H{o}s-R\'enyi random graphs via tree counting}
\author{ \textbf{Ghurumuruhan Ganesan}
\thanks{E-Mail: \texttt{gganesan82@gmail.com} } \\
\ \\
New York University, Abu Dhabi}
\date{}
\maketitle

\begin{abstract}
Consider the complete graph \(K_n\) on \(n\) vertices where each edge \(e\) is independently open with probability \(p_n(e)\) or closed otherwise. Here \(\frac{C-\alpha_n}{n} \leq p_n(e) \leq \frac{C+\alpha_n}{n}\) where \(C > 0\) is a constant not depending on~\(n\) or~\(e\) and \(0 \leq \alpha_n  \longrightarrow 0\) as \(n \rightarrow \infty.\) The resulting random graph~\(G\) is inhomogenous and we use a tree counting argument to establish phase transition in \(G.\) We also obtain that the critical value for phase transition is one in the following sense. For \(C < 1,\) all components of \(G\) are small (i.e. contain at most \(M\log{n}\) vertices) with high probability, i.e., with probability converging to one as \(n \rightarrow \infty.\) For \(C > 1,\) with high probability, there is at least one giant component (containing at least \(\epsilon n\) vertices for some \(\epsilon > 0\)) and every component is either small or giant. For \(C > 8,\) with positive probability, the giant component is unique and every other component is small. As a  consequence of our method, we directly obtain the fraction of vertices present in the giant component in the form of an infinite series.

\vspace{0.1in} \noindent \textbf{Key words:} Erd\H{o}s-R\'enyi random graphs, phase transition, tree counting argument.

\vspace{0.1in} \noindent \textbf{AMS 2000 Subject Classification:} Primary:
60J10, 60K35; Secondary: 60C05, 62E10, 90B15, 91D30.
\end{abstract}

\bigskip

\setcounter{equation}{0}
\renewcommand\theequation{\thesection.\arabic{equation}}
\section{Introduction} \label{intro}
Consider the complete graph \(K_n\) on \(n\) vertices \(\{1,2,...,n\}\) and corresponding edges \(e_1,\ldots,e_m\) where \(m = {n \choose 2}.\) An edge~\(e_i\) is independently open with probability \(p_n(e_i),\) and closed otherwise. Throughout we assume that
\begin{equation}\label{edge_prob}
\frac{C-\alpha_n}{n} \leq \min_{1 \leq i \leq m} p_n(e_i) \leq \max_{1 \leq i \leq m} p_n(e_i) \leq \frac{C+\alpha_n}{n}
\end{equation}
where \(0 \leq \alpha_n \longrightarrow 0\) as \(n \rightarrow \infty.\) The resulting random graph \(G\) is an Erd\H{o}s-R\'enyi (ER) graph (Janson et al (2000)).  Strictly speaking \(G\) is one version of an ER random graph and in the original paper, Erd\H{o}s and R\'enyi (1960) have studied phase transition in a slightly different model of random graphs.

When \(p_n(e_i) = \frac{C}{n}\) for all \(i,\) the probabilities for an edge begin open are the same and the resulting random graph is homogenous. Phase transitions in homogenous graphs have been studied in great detail and the common approaches include a combination of vertex exploration, random walks and branching processes (see Janson~et~al~(2000), Durrett (2007), Alon and Spencer (2008) and references therein) and usually look at stochastic domination of the exploration process from above and below by branching processes. For a more recent comprehensive account, we refer to the monographs of Hofstad~(2016). For example, the results related to phase transitions are treated in Chapter \(4,\) Theorems \(4.4\)-\(4.8\) in Hofstad~(2016).

When the edge probabilities are not the same, the resulting random graph~\(G\) is inhomogenous and the above analysis is not directly applicable. We illustrate with an example in Section~\ref{example}. To study phase transitions in inhomogenous random graphs, we use in this paper a combinatorial tree counting argument to estimate the sizes of small and mid size components. To make the proof self contained, we give the proofs for both subcritical and super critical regimes. 

Tree counting arguments have been used before in the analysis of random graphs. For example, Bollobas (2001) (Chapter~\(7\)), Durrett (2007) (Chapter~\(2\)) and Alon and Spencer (2008) (Chapter~\(11\)) have all used tree counting arguments in various forms for different purposes. Bollobas (2001) uses the tree counting with coarser estimates to analyse the connectivity regime of homogenous random graphs where each edge is independently open with probability \(p = \frac{C\log{n}}{n}\) for some constant \(C > 0.\) Durrett (2007) obtains estimates on the number of tree components of a fixed size using the tree counting argument for \(p = \frac{C}{n}.\) Alon and Spencer (2008) use the tree counting argument with \(p = \frac{C}{n}\) for comparison with a Poisson branching process.

Our argument is different from the above in the sense that we obtain more precise tree counting estimates for the small component sizes when edges are open with probability of the order of \(\frac{1}{n}\) (see Lemmas~\ref{c1_est1_lem} and~\ref{cross}, Section~\ref{prelim}). We use the above estimates to indirectly deduce the presence of a giant component. An important consequence of our method is that for the supercritical regime of \(C > 1,\) we obtain the analytical solution for the fraction~\(q(C)\) of vertices present in giant components in the form of an infinite series. For details, we refer to the proof of Theorem~\ref{th1}, Section~\ref{pf_sup}.

\subsubsection*{Model Description}
We briefly describe the probability space first. For \(1 \leq i \leq m,\) we define the state of the edge~\(e_i \in K_n\) on the probability space \((\Omega_i, {\cal F}_i, \mathbb{P}_i)\) where \(\Omega_i = \{0,1\}, {\cal F}_i = \{\emptyset, \{0\}, \{1\}, \{0,1\}\}\) and \(\mathbb{P}_i\) denotes the Bernoulli probability measure with \(\mathbb{P}_i(\{1\}) = p_n(e_i) =  1-\mathbb{P}_i(\{0\}).\) Thus open edges are assigned a value of \(1\) and closed edges are assigned a value of \(0.\) We define the random graph \(G\) on the probability space \((\Omega, {\cal F}, \mathbb{P})\) where \(\Omega = \otimes_{i=1}^{m} \Omega_i, {\cal F} = \otimes_{i=1}^{m}{\cal F}_i\) and \(\mathbb{P} = \otimes_{i=1}^{m} \mathbb{P}_i.\)

Throughout the paper, we study open components of the graph \(G\) and we give a brief description. Let \(e_{i}\) and \(e_j, j \neq i\) be two edges in \(K_n.\) We say~\(e_i\) and~\(e_j\) are \emph{adjacent} if they share an endvertex. We say that a sequence of distinct edges \(P = (f_1,f_2,\ldots,f_k), f_i \subset \{e_j\}\) is a \emph{path} in \(K_n\) if the edge~\(f_i\) is adjacent to the edge~\(f_{i+1}\) for every \(1 \leq i \leq k-1.\) We say that~\(P\) is an \emph{open path} (in \(G\)) if \(P\) is a path and every~\(f_i, 1 \leq i \leq n,\) is open. Let \(v_1\) be the endvertex of the edge~\(f_1\) not common to~\(f_2\) and let \(v_k\) be the endvertex of \(f_k\) not common to~\(f_{k-1}.\) We say that~\(v_1\) and~\(v_k\) are \emph{endvertices} of the path~\(P.\)

Let \(1 \leq i\neq j \leq n\) be fixed. We say that vertex \(i\) is connected to vertex \(j\) by an open path if there is an open path \(P_{ij}\) containing \(i\) and \(j\) as endvertices. Let~\({\cal E}_i\) be the set of all vertices \(v, 1 \leq v \leq n,\) connected to the vertex~\(i\) by an open path. We assume that \(i \in {\cal E}_i\) and define \({\cal E}_i\) to be the \emph{open component} of the graph~\(G\) containing the vertex~\(i.\) We also refer to open components simply as components.

For \(C > 0,\) let
\begin{equation}\label{del_def2}
\delta  = \delta(C) := C- 1 - \log{C}.
\end{equation}
We have that \(\delta(1) = 0\) and \(\delta(C) >0\) for \(C \neq 1.\) We are interested in the size of components of \(G\) for the case \(C < 1\) (subcritical) and the case \(C > 1\) (supercritical) which we discuss separately below.

\subsection*{Subcritical case}
We have the following result for \(C< 1.\) For \(1 \leq i \leq n,\) let \({\cal E}_i\) denote the open component of the graph \(G\) containing the vertex~\(i\) defined in the previous subsection. For integer \(M \geq 1\) let
\begin{equation}\label{h_2_def}
H_1 = H_1(M) := \bigcap_{1 \leq j \leq n}\{\#{\cal E}_j \leq M \log{n}\}
\end{equation}
denote the event that the size of every component is at most \(M \log{n}.\) Throughout, the size of a component refers to the number of vertices present in the component.

\begin{Theorem}\label{th2} Fix \(C < 1\) and let \(\delta = \delta(C) > 0\) be as in (\ref{del_def2}). Fix \(M > \frac{1}{\delta}\)
and let \(\theta < \min(M\delta - 1, 1)\) be fixed.
There is a positive constant \(L_1 = L_1(C, M, \theta)\) so that
\begin{equation}\label{h_2_est}
\mathbb{P}(H_1(M)) \geq 1- \frac{1}{n^{\theta}}
\end{equation}
for all \(n \geq L_1.\)
\end{Theorem}
Thus with high probability, (i.e., with probability converging to one as \(n \rightarrow \infty\)), every component has size at most a constant multiple of~\(\log{n}.\)


For values of \(C\) slightly farther away from the critical value of one, we have stronger rates of decay.
\begin{Theorem}\label{th112} Fix \(C < e^{-1}\) and \(\theta > 0.\) There are positive constants \(M = M(C,\theta)\) and \(L_2 = L_2(C,\theta,M)\) so that
\begin{equation}\label{e_m_111}
\mathbb{P}(H_1(M)) \geq 1-\frac{1}{n^{\theta}}
\end{equation}
for all \(n \geq L_2.\)
\end{Theorem}

\subsection*{Supercritical case}
To study the supercritical case of \(C > 1,\) we first have some preliminary definitions. For integer \(r \geq 1,\) let \(T_r\) denote the number of labelled trees on~\(r\) vertices.
We set \(T_1 = 1\) and for \(r \geq 2,\) we recall the Cayley formula (see van Lint and Wilson (2010)) for \(T_r\) as \(T_r = r^{r-2}\) for \(r \geq 2.\)
In this paper, we do not require the use of the exact Cayley formula and therefore, we simply treat \(T_r\) as a sequence.

For \(C > 0,\) define
\begin{equation}\label{q_0_def}
q = q(C) := \sum_{r=1}^{\infty} \frac{T_re^{-r}}{C(r-1)!}e^{-\delta r} = \sum_{r=1}^{\infty} \frac{T_rC^{r-1}e^{-Cr}}{(r-1)!}
\end{equation}
where \(\delta = \delta(C) \geq 0\) is as in (\ref{del_def2}). For the supercritical case, the function~\(q(C)\) plays a crucial
role in determining the fraction of vertices in the giant component of \(G.\)

The following result collects together the important properties of \(q(C)\) needed for the proofs.
\begin{Theorem}\label{th_q0}The function \(q(C)\) satisfies the following properties:\\
\((i)\) For \(C > 0,\) we have
\begin{equation}\label{q0_up_bd123}
0 < q(C) \leq 1.
\end{equation}
\((ii)\) If \(C \leq 1,\) then
\begin{equation}\label{q0_eq1}
q(C) = 1.
\end{equation}
\((iii)\) For \(C > 1,\) the function \(q(C)\) is strictly decreasing in \(C\) and
\begin{equation}\label{q0_eq2}
0 < q(C) < 1.
\end{equation}
\((iv)\)  For \(C > 0,\) we have that \(q = q(C)\) satisfies the relation
\begin{equation}\label{extinct}
q = e^{-C(1-q)}.
\end{equation}
\end{Theorem}
Using a combinatorial approach, like for e.g. generating functions, we obtain the property \((iv)\) that~\(q = q(C)\) satisfies the relation~(\ref{extinct}) (see e.g. van Lint and Wilson~(2001), Grimmett (1980)). The term \(q\) satisfying~(\ref{extinct}) is also the extinction probability of a Poisson
branching process with mean offspring size \(C\) (see Durrett (2007)). Thus the expansion \(q(C)\) in (\ref{q_0_def}) is the analytical solution for the extinction probability. We give a probabilistic proof of properties \((i)-(iii)\) in Section~\ref{pf_q0} using the properties of random graphs. 

Fix \(C > 1\) and \(\epsilon > 0.\) For \(i \geq 1,\) let
\begin{equation}\label{y_i_def}
Y_i := \ind(\#{\cal E}_i > \epsilon n)
\end{equation}
denote the indicator function of the event that vertex~\(i\) belongs to a giant component, i.e., a component containing more than~\(\epsilon n\) vertices. We then have
\begin{equation}\label{y_def}
Y_n(\epsilon) := \sum_{i=1}^{n} Y_i
\end{equation}
denotes the sum of sizes of all giant components. We use the term giant component to roughly mean any component which contains at least a constant fraction \(f \in (0,1)\) of the~\(n\) vertices. We make the notion precise based on the context.

For \(0 < \gamma < 1\) and define the event
\begin{equation}\label{v_i_def}
V(\gamma,\epsilon) := \{(1-q(C)-\gamma)n \leq Y_n(\epsilon) \leq (1-q(C)+\gamma)n\}.
\end{equation}
For \(M > 0\) let
\begin{equation}\label{h_1_def}
H_2 = H_2(M,\gamma,\epsilon) := V(\gamma,\epsilon) \bigcap \bigcap_{1 \leq j \leq n}\left(\{\#{\cal E}_j \leq M \log{n}\}\bigcup \{\#{\cal E}_j > \epsilon n\} \right)
\end{equation}
denote the event that \(V(\gamma,\epsilon)\) occurs and every component is either giant or small. Here we say that a component is small if its size is at most~\(M\log{n}.\) We have the following result.
\begin{Theorem}\label{th1} Fix \(C > 1\) and let \(\delta  = \delta(C) > 0\) be as in (\ref{del_def2}). Fix \(M > \frac{1}{\delta}\) and \(0 < \gamma  < 1.\) There exists~\(\epsilon_1 = \epsilon_1(C,M,\gamma) > 0\) such that the following holds for all
\(0 < \epsilon < \epsilon_1(C).\) There is a positive constant \(L_3 = L_3(C,M,\gamma,\epsilon)\) so that
\begin{equation}\label{h_1_est}
\mathbb{P}(H_2) \geq 1-\gamma
\end{equation}
for all \(n \geq L_3.\)
\end{Theorem}
Since \(q(C) < 1\) for \(C > 1\) and \(0 < \gamma < 1\) is arbitrary, we have with arbitrarily large probability that there exists a giant component of~\(G.\)

\subsection*{Uniqueness of the giant component}
In Theorem~\ref{th1}, we have established that for any \(C > 1,\) with high probability there is at least one giant component, i.e., a component whose size is larger than~\(\epsilon n\) for some constant \(\epsilon > 0.\) We now see that for values of \(C\) slightly farther away from one, the giant component is unique with positive probability.

For \(0 < \gamma <1\) and for \(1 \leq i \leq n,\) let
\begin{equation}\label{v_i_def2}
V_i = V_i(\gamma) := \{(1-q(C)-\gamma)n \leq \#{\cal E}_i \leq (1-q(C)+\gamma)n\}
\end{equation}
denote the event that the size of the component \({\cal E}_i\) containing the vertex \(i\) is in the range \([(1-q(C)-\gamma)n,(1-q(C)+\gamma)n].\)
Here \(q(C)\) is as in (\ref{q_0_def}).
For \(M > 0\) define
\begin{equation}\label{h_11_def}
H_3 = H_3(M,\gamma) := \bigcup_{i=1}^{n} V_i \bigcap \bigcap_{1 \leq j \leq n, j \neq i}\{\#{\cal E}_j \leq M \log{n}\}
\end{equation}
to be event that there exists a unique giant component with size in the range \([(1-q(C)-\gamma)n, (1-q(C)+\gamma)n]\) and the size of every other component is at most~\(M \log{n}.\) For positive numbers \(C, \epsilon\) and \(\omega,\) let \(\delta_0(C,\epsilon,\omega) = \delta(C) - C\epsilon-\omega.\) We also use the definition of \(\delta_0\) in future results.
\begin{Theorem}\label{th111} There is a unique \(C_0 > 1\) such that \(\delta_0\left(C_0 ,\frac{1}{2},0\right) = 0.\) We have that \(C_0  < 8.\)

Fix \(C > C_0\) and let \(\delta = \delta(C)> 0\) be as in (\ref{del_def2}). Fix \(M > \frac{1}{\delta}\) and let \(\alpha > 0\) be such that \(\delta_0\left(C,\alpha + \frac{1}{2},0\right) > 0.\) Fix \(0 < \gamma <1.\) There is a positive constant \(L_4 = L_4(C,\alpha,\gamma, M)\) so that
\begin{equation}\label{e_m_111}
\mathbb{P}(H_3(M,\gamma)) \geq 1-e^{-\alpha C}-\gamma
\end{equation}
for all \(n \geq L_4.\)
\end{Theorem}
For any \(C  > 8,\) we therefore have with positive probability that there exists a unique giant component.






The paper is organized as follows: In Section~\ref{example}, we provide an example of an inhomogenous random graph. The three following sections obtain preliminary estimates. In Section~\ref{prelim}, we collect the tree counting estimates for non giant components, i.e., components with size at most \(\epsilon n\) for some constant~\(\epsilon > 0.\) In Section~\ref{mid_comp}, we obtain estimates on the mid size components whose size lies in the range \([M\log{n}+1, \epsilon n]\) for some constant \( M >0.\) We then obtain mean and variances estimates for \(Z_n(\epsilon) = n - Y_n(\epsilon)\) (see (\ref{y_def})) in Section~\ref{mean_zn} needed for the proofs of the main Theorems.

Using the preliminary estimates of the above Sections, we first prove Theorem~\ref{th2} regarding the subcritical case \(C<1\) in Section~\ref{pf_sub}. We then provide a probabilistic proof for properties \((i)-(iii)\) of Theorem~\ref{th_q0} in Section~\ref{pf_q0}. In Section~\ref{pf_sup}, we prove the supercritical case \(C > 1\) of Theorem~\ref{th1} and finally, in Section~\ref{pf11}, we prove Theorems~\ref{th111} and~\ref{th112}.

\setcounter{equation}{0}
\renewcommand\theequation{\thesection.\arabic{equation}}
\section{Inhomogenous random graph example}\label{example}
Suppose that the edge probabilities \(p_n(e) \in \left\{\frac{C-\alpha_n}{n}, \frac{C}{n}\right\}\) for all edges \(e \in K_n\) and suppose that \[\#\left\{e \in K_n : p_n(e) = \frac{C-\alpha_n}{n}\right\} = \frac{1}{8}n^2.\] We recall that \(K_n\) is the complete graph on \(n\) vertices. Thus there are \(\frac{1}{8}n^2\) edges in~\(K_n\) with \(p_n(e) = \frac{C-\alpha_n}{n}\) and the rest have \(p_n(e) = \frac{C}{n}.\)

To compare homogenous and inhomogenous random graphs, we perform a coupling procedure as follows. Let \(\{X(e)\}_{e \in K_n}\) be a set of independent and identically distributed (i.i.d) random variables uniformly distributed in \([0,1]\) and defined on the probability space \((\Xi,{\cal G}, \mathbb{P}_c).\) For a particular realization of \(\{X(e)\}_{e \in K_n},\) define the graphs \[G^{-} = \left\{e \in K_n : X(e) \leq \frac{C-\alpha_n}{n}\right\}\] and \[G = \left\{e \in K_n : X(e) \leq p_n(e)\right\}.\] Thus~\(G\) is the inhomogenous random graph where edge \(e\) is open with probability~\(p_n(e)\) independently of all the other edges. The graph~\(G^{-}\) is the homogenous random graph obtained when every edge is independently open with probability \( \frac{C-\alpha_n}{n}.\)

The phase transition for the homogenous graph~\(G^{-}\) essentially involves the study of two properties (Durrett (2007)). \((a)\) The existence of a unique giant component containing roughly~\((1-q(C))n\) vertices and~\((b)\)~No middle ground in the sense that every component other than the giant component has at most \(D\log{n}\) vertices for some constant \(D > 0.\) As described in Section~\ref{intro}, the term~\(q(C)\) is the probability of extinction for the Poisson branching process.

By construction, the graph \(G^{-} \subseteq G.\) However, we argue below that even if the homogenous graph \(G^{-}\) satisfies the phase transition properties~\((a)\) and~\((b)\) above, the inhomogenous random graph \(G\) need not. We use the following two estimates.\\
\((a1)\) Let \[{E}_{dif} = \left\{e \in K_n : \frac{C-\alpha_n}{n} < X_n(e) \leq p_n(e)\right\}\] be the set of edges which are open in \(G\) but closed in \(G^{-}.\) Let \(R_{dif}\) be the number of edges in \(E_{dif}.\) If \(\alpha_n = \frac{16}{\sqrt{n}},\) then we have that \[\mathbb{P}_c\left(R_{dif} \in [\sqrt{n},3\sqrt{n}]\right) \geq 1-\frac{2}{\sqrt{n}}\] for all \(n \geq 2.\)\\
\((a2)\) For \(1 \leq i \leq n,\) let \(d_i\) denote the degree of vertex~\(i\) in the graph \(G.\) We have that
\[\mathbb{P}_c\left(\sup_{1 \leq i \leq n} d_i  \leq 3\log{n}\right) \geq 1-\frac{D}{n^2}\] for some constant \(D >0\) and for all \(n \geq 2.\)\\

From \((a1)-(a2)\) we have that with high probability (i.e. with probability converging to one as \(n \rightarrow \infty\)), the random variable \(R_{dif}\) lies in the range~\([\sqrt{n},3\sqrt{n}]\) and the maximum degree of a vertex in the graph \(G\) is at most \(3\log{n}.\)\\
\emph{Proof of \((a1)-(a2)\)}: For \((a1)\) we proceed as follows. The random variable~\(\#E_{dif}\) is binomially distributed with parameters \(m = \frac{1}{8} n^2\) and \(p_b = \frac{\alpha_n}{n}.\) For \(\alpha_n = \frac{16}{\sqrt{n}},\) we have that the mean \[\mathbb{E}_c(R_{dif}) = mp_b = \frac{\alpha_n}{n} \frac{1}{8}n^2 = 2\sqrt{n}\] and the variance \(var(R_{dif})  = mp_b(1-p_b) \leq mp_b = 2\sqrt{n}.\) Here \(\mathbb{E}_c\) denotes the expectation operator corresponding to the probability measure \(\mathbb{P}_c.\) Using the Chebychev inequality, we therefore have that
\begin{equation}\label{r_dif_est}
\mathbb{P}_c\left(|R_{dif} - \mathbb{E}R_{dif}| > \sqrt{n}\right) \leq \frac{var(R_{dif})}{n} \leq \frac{2\sqrt{n}}{n} = \frac{2}{\sqrt{n}}
\end{equation}
This proves \((a1).\)

We prove \((a2)\) as follows. Suppose \(d_i = \sum_{j=1}^{n} Z_{i,j}\) denotes the degree of vertex~\(i,\) where \(Z_{i,j} = 1\) if vertices \(i\) and \(j\) are joined by an edge and zero otherwise. The random variables \(\{Z_{i,j}\}_{1 \leq j \leq n}\) are independent and for \(1 \leq j \leq n,\) we have \[\mathbb{E}_cZ_{i,j} =p_{i,j} \leq p_u = \frac{C+\alpha_n}{n} \leq \frac{C+1}{n}\] for all \(n \geq n_0\) large, using (\ref{edge_prob}). The number \(n_0\) does not depend on the choice of~\(i\) or~\(j.\) For \(s > 0\) we therefore have
\[\mathbb{E}_ce^{sZ_{i,j}} = 1-p_{i,j} + e^{s}p_{i,j}  = 1+(e^s-1)p_{i,j} \leq 1+(e^s-1)p_u\] and \[\mathbb{E}_ce^{sd_i} = \prod_{j=1}^{n} \mathbb{E}_ce^{sZ_{i,j}} \leq (1+(e^s-1)p_u)^{n}\leq \exp\left((e^{s}-1)(C+1)\right)\] where the last estimate is obtained using \(1+x \leq e^{x}\) for \(x > 0.\)

Setting \(s = 1\) and using the Markov inequality we have that
\begin{equation}\nonumber
\mathbb{P}_c\left(d_i \geq 3\log{n}\right) \leq \frac{\mathbb{E}e^{d_i}}{n^3} \leq  \frac{D}{n^3}
\end{equation}
for all \(n \geq n_0\) and for \(D = \exp((e-1)(C+1)).\) Thus
\begin{equation}\label{di_est}
\mathbb{P}_c\left(\bigcup_{1 \leq i \leq n} \{d_i \geq 3\log{n}\}\right) \leq \frac{D}{n^2}
\end{equation}
for all \(n \geq n_0\) and so with high probability, every vertex in the graph \(G\) has degree at most \(3\log{n}.\) This proves \((a2).\) \(\qed\)

We use properties \((a1)-(a2)\) to compare the inhomogenous random graph~\(G\) with the homogenous random graph~\(G^-.\) Suppose that the term \(R_{dif} \in [\sqrt{n},3\sqrt{n}]\) and every vertex in the graph \(G\) has degree at most \(3\log{n}.\) By construction, we have that \(G^{-} \subset G\) and the random variable~\(R_{dif}\) is the number of edges which are open in \(G\) but closed in \(G^-.\) Even if the homogenous graph~\(G^{-}\) satisfies the phase transition properties~\((a)\) and~\((b)\) above, the graph~\(G\) could have the following component structure. The giant component in \(G \supset G^-\) is the same giant component as in~\(G^-.\) There is a component~\(C_{mid}\) consisting of~\(x \in [\sqrt{n},3\sqrt{n}]\) edges all of which belong to \({\cal E}_{dif},\) i.e., the set of edges open in~\(G\) but closed in~\(G^{-}\) and the rest of the components of \(G\) have size at most \(D \log{n}.\)

Let \(\#C_{mid}\) be the size, i.e., the number of vertices in the component \(C_{mid}.\) We have that
\begin{equation}\label{c_mid_size}
\frac{2\sqrt{n}}{3\log{n}} \leq \#C_{mid} \leq 6\sqrt{n}
\end{equation}
and so that graph \(G\) does not satisfy the phase transition properties~\((a)-(b).\) \\
\emph{Proof of (\ref{c_mid_size})}:
For the lower bound, we use the fact that \(\sum_{v \in C_{mid}} d_v = 2x  \geq 2\sqrt{n}\) where \(d_v\) denotes the degree of vertex \(v \in C_{mid}.\) Since the degree of every vertex in~\(G\) is at most~\(3\log{n},\) we also have that \(\sum_{v \in C_{mid}} d_v \leq 3\log{n} \#C_{mid}.\)  Thus \(\#C_{mid} \geq \frac{2\sqrt{n}}{3\log{n}}.\)

For the upper bound, we use the fact that there are \(x \leq 3\sqrt{n}\) edges in~\(C_{mid}\) and so the number of vertices in \(C_{mid}\) is at most \(6\sqrt{n}.\) \(\qed\)

An analogous argument as above holds if we use the homogenous random graph \(G^+ \supset G\) obtained by allowing every edge to be independently open with probability~\(\frac{C}{n}.\)

\setcounter{equation}{0}
\renewcommand\theequation{\thesection.\arabic{equation}}
\section{Tree counting estimates}\label{prelim}
For \(C > 0\) let \(\delta = \delta(C)\) be as in (\ref{del_def2}).
For \(\omega  >0, \epsilon \in (0,1)\) and \(i = 0,1,2,\) define \(\delta_i = \delta_i(C,\epsilon,\omega)\) as
\begin{equation}\label{delta0_def}
\delta_0= \delta - C\epsilon - \omega,\delta_1 = \delta - 2C\epsilon - \omega \text{ and } \delta_2 =\delta - \log(1-\epsilon) + \omega.
\end{equation}
For any \(\omega > 0, \epsilon \in (0,1),\) we have that \(\delta_2(C,\epsilon,\omega) > 0.\) 



\subsection*{Non giant components estimate}
For \(i \geq 1,\) let \({\cal E}_i\) denote the component containing the vertex \(i.\)  To study non-trivial components of \({\cal E}_i,\) we let \(T_1 = 1\) and as before, for \(r \geq 2\) we let~\(T_r\) be the number of labelled trees on~\(r\) vertices. We have the following Lemma that obtains upper and lower bounds concerning the size of non giant components, i.e., components whose size is at most \(\epsilon n\) for some constant \(\epsilon  >0.\)

\begin{Lemma}\label{c1_est1_lem}Fix \(C \neq 1\)  and fix \(\epsilon,\omega > 0\) so that both \(\delta_0\) and \(\delta_2\) defined in~(\ref{delta0_def}) are positive. There is a positive constant \(N = N(C,\epsilon,\omega)\) such that the following two statements hold. For all \(n \geq N, 1 \leq i \leq n\) and \(1 \leq r \leq \epsilon n,\) we have
\begin{eqnarray}
\mathbb{P}\left(\#{\cal E}_i = r\right) \leq \frac{T_r e^{-r}}{C(r-1)!}e^{-\delta_0 r}.\label{c1_est1}
\end{eqnarray}
For all \(n \geq N, 1 \leq i \leq n\)  and \(1 \leq r \leq \epsilon n,\) we have
\begin{eqnarray}
\mathbb{P}\left(\#{\cal E}_i = r\right) \geq \frac{1}{C(1-\epsilon)}\frac{T_r e^{-r}}{(r-1)!}e^{-\delta_2 r}e^{-2\epsilon + \frac{2\omega}{3}}.\label{c1_est2}
\end{eqnarray}
\end{Lemma}
We need the lower bound in (\ref{c1_est2}) to estimate the fraction of  nodes present in the giant component for \(C > 1.\)

\emph{Proof of Lemma~\ref{c1_est1_lem}}: Let \(C_u = C+\alpha_n, C_d = C-\alpha_n\) and \(p_u = \frac{C_u}{n}\) and \(p_d = \frac{C_d}{n}.\)

We prove the upper bound in (\ref{c1_est1}) first. If the vertex \(i\) is isolated, then every edge containing~\(i\) as an endvertex is closed. Since every edge is closed with probability at most \(1-p_d\) and  at least \(1-p_u,\) we have
\begin{equation}
(1-p_u)^{n-1} \leq \mathbb{P}(\#{\cal E}_i = 1) \leq (1-p_d)^{n-1}. \label{eq_c_11}
\end{equation}
For components with larger size, we argue for \(i = 1\) and then generalize for all~\(i.\) Suppose now that \(i = 1\) and the component~\({\cal E}_1\) contains \(r \geq 2\) vertices. There is a random tree \({\cal J}\) contained in~\({\cal E}_1\) with the same vertex set as~\({\cal E}_1\) and containing~\(r-1\) edges, each of which is open. Moreover, every edge with one endvertex belonging to \({\cal J}\) and the other endvertex not in \({\cal J},\) is closed. The number of such edges is~\(r(n-r).\) Let \({\cal T}_r\) denote the set of all trees with vertex set \(\{1,2,\ldots,r\}.\) For a fixed tree \({\cal T} \in {\cal T}_r,\) we therefore have
\begin{equation}
\mathbb{P}\left(\{{\cal E}_1 = \{1,2,\ldots,r\}\} \cap \{{\cal T} \subseteq {\cal E}_1\}\right) \leq p_u^{r-1}(1-p_d)^{r(n-r)}. \label{eq_c_bg21}
\end{equation}
The final expression is because every edge in the graph \(K_n\) is open with probability at most \(p_u\) and closed with probability at most \(1-p_d\)~(see (\ref{edge_prob})). Summing over all possible choices of \({\cal T},\) we have
\begin{equation}
\mathbb{P}\left(\{{\cal E}_1 = \{1,2,\ldots,r\}\}\right) \leq T_r p_u^{r-1}(1-p_d)^{r(n-r)} \label{eq_c_bg22a}
\end{equation}
where as before \(T_r\) denotes the number of labelled trees on \(r\) vertices.

The estimate (\ref{eq_c_bg22a}) is for a particular choice of vertex set for the component \({\cal E}_1.\) The total number of choices for the vertex set of~\({\cal E}_1\) is the number of ways of choosing \(r-1\) vertices (apart from the vertex \(1\)) out of the remaining \(n-1\) vertices. We therefore have from~(\ref{eq_c_bg22a}) that
\begin{equation}
\mathbb{P}\left(\#{\cal E}_1 = r\right) \leq  {n-1 \choose r-1} T_r p_u^{r-1}(1-p_d)^{r(n-r)}. \nonumber
\end{equation}
The above argument holds for all \(1 \leq i \leq n\) and so
\begin{equation}
\mathbb{P}\left(\#{\cal E}_i = r\right) \leq  {n-1 \choose r-1} T_r p_u^{r-1}(1-p_d)^{r(n-r)}. \label{eq_c_bg2}
\end{equation}

Using \(p_u = \frac{C_u}{n}\) we evaluate the product of the first and third terms as
\begin{eqnarray}
{n-1 \choose r-1} p_u^{r-1} = \frac{(n-1)!}{(r-1)!(n-r)!}\frac{C_u^{r-1}}{n^{r-1}} \leq \frac{C_u^{r-1}}{(r-1)!}.\label{eq1}
\end{eqnarray}
The final estimate holds since \(\frac{(n-1)!}{(n-r)!n^{r-1}} = \frac{(n-1)\ldots(n-r+1)}{n^{r-1}}  \leq 1.\)

To evaluate the last term in (\ref{eq_c_bg2}), we again use \(1-x \leq e^{-x}\) for \(x >0\) to get
\begin{equation}\label{temp1}
(1-p_d)^{r(n-r)} \leq \exp\left(-\frac{C_dr(n-r)}{n}\right) = e^{-C_dr}\exp\left(\frac{C_dr^2}{n}\right).
\end{equation}
Substituting (\ref{eq1}) and (\ref{temp1})
into (\ref{eq_c_bg2}) we get
\begin{eqnarray}
\mathbb{P}\left(\#{\cal E}_i = r\right) \leq T_r\frac{C_u^{r-1}}{(r-1)!}e^{-C_dr}\exp\left(\frac{C_dr^2}{n}\right).\label{c1_r_raw}
\end{eqnarray}

For \(r \leq \epsilon n,\) we have \(\frac{C_dr^2}{n} \leq \epsilon C_d r\) and
therefore
\begin{eqnarray}
\mathbb{P}\left(\#{\cal E}_i = r\right) \leq T_r\frac{C_u^{r-1}e^{-C_dr}}{(r-1)!}e^{C_d\epsilon r}\label{est_c1_ra}
\end{eqnarray}
Using \(C_u = C+\alpha_n\) and \(C_d = C-\alpha_n,\) we have that
\begin{eqnarray}
C_u^{r-1} e^{-C_d r}e^{C_d\epsilon r} &=& C^{r-1}(1+\alpha_n C^{-1})^{r-1} e^{-Cr} e^{\alpha_n r} e^{C\epsilon r} e^{-\epsilon \alpha_n r} \nonumber\\
&\leq& C^{r-1}e^{\alpha_nC^{-1}(r-1)} e^{-Cr} e^{\alpha_n r} e^{C\epsilon r} e^{-\epsilon \alpha_n r} \nonumber\\
&\leq& C^{r-1}e^{-Cr} e^{C\epsilon r} e^{\alpha_nC^{-1}(r-1)} e^{\alpha_n r}. \label{temp_cu_cd}
\end{eqnarray}
The middle inequality is because \(1+x \leq e^{x}\) for all \(x > 0.\)

Since \(r \geq 1,\) we have that \(\alpha_nC^{-1}(r-1) \leq \alpha_nC^{-1} r\) and so we have
\begin{equation}\label{omega_est}
\alpha_nC^{-1} (r-1) + \alpha_n r \leq (C^{-1}+1)\alpha_n r \leq \omega r
\end{equation}
for all \(n\) large.  Here \(\omega > 0\) is as in the statement of the Lemma and the final estimate above holds since \(\alpha_n \longrightarrow 0\) as \(n \rightarrow \infty\) and so \((C^{-1}+1)\alpha_n < \omega\) for all \(n\) large.


Substituting (\ref{omega_est}) into (\ref{temp_cu_cd}), we have that
\begin{equation}
C_u^{r-1} e^{-C_d r}e^{C_d\epsilon r} \leq C^{r-1}e^{-Cr} e^{C\epsilon r} e^{\omega r} = \frac{1}{C}e^{-\delta_0 r}e^{-r} \label{cross_later}
\end{equation}
for all \(n\) large. Here \(\delta_0 > 0\) is as in (\ref{delta0_def}).  Substituting~(\ref{cross_later}) into~(\ref{est_c1_ra}) we get (\ref{c1_est1}).

\emph{Proof of (\ref{c1_est2})}: The proof is analogous as above. As before, we argue for \( i =1\) and then generalize for all~\(i.\) Let \({\cal T}_r\) denote the set of labelled trees with vertex set \(\{1,2,\ldots,r\}.\) We then have
that
\begin{equation}\label{t_est1}
\mathbb{P}({\cal E}_1 = \{1,2,\ldots,r\}) \geq \sum_{{\cal T} \in {\cal T}_r} \mathbb{P}({\cal E}_1 = {\cal T}).
\end{equation}
For any fixed tree \({\cal T} \in {\cal T}_r,\) suppose \(\{{\cal E}_1 = {\cal T}\}\) occurs. We then have that the \(r-1\) edges of \({\cal T}\) are open and every edge with one vertex in \({\cal T}\) and one vertex outside is closed. The number of such edges is \(r(n-r).\) Moreover, since \({\cal E}_1\) has exactly \(r-1\) edges, the remaining \({r \choose 2} - r+1\) edges in \(K_r\) are closed. Here \(K_r\) is the complete graph with vertex set \(\{1,2,\ldots,r\}.\) Using the fact that edges are open with probability at least \(p_d\) and closed with probability at least \(1-p_u,\) we have that
\begin{equation}\label{t_est2}
\mathbb{P}({\cal E}_1 = {\cal T}) \geq p_d^{r-1}(1-p_u)^{r(n-r) + {r \choose 2} - r+1}
\end{equation}
and since the probability is the same for any tree in \({\cal T}_r,\) we substitute (\ref{t_est2}) into (\ref{t_est1}) to get
\begin{equation}\label{t_est3}
\mathbb{P}({\cal E}_1 = \{1,2,\ldots,r\}) \geq T_r p_d^{r-1}(1-p_u)^{r(n-r) + {r \choose 2} - r+1}
\end{equation}
where \(T_r\) is the number of labelled trees on \(r\) vertices. 


Since there are \({n-1 \choose r-1}\) ways to choose the remaining \(r-1\) vertices for \({\cal E}_1,\) we therefore have
\begin{equation}\nonumber
\mathbb{P}(\#{\cal E}_1 = r) \geq {n-1 \choose r-1} T_r p_d^{r-1}(1-p_u)^{r(n-r) + {r \choose 2} - r+1}.
\end{equation}
Again setting \(T_1 = 1\) and using (\ref{eq_c_11}), the above bound also holds for \(r = 1.\) 
The above argument holds for all \(i\) and so we have
\begin{equation}\label{t_est4}
\mathbb{P}(\#{\cal E}_i = r) \geq {n-1 \choose r-1} T_r p_d^{r-1}(1-p_u)^{r(n-r) + {r \choose 2} - r+1}.
\end{equation}

It remains to simplify (\ref{t_est4}) to get (\ref{c1_est2}).
Using \(p_d = \frac{C_d}{n},\) we have
\begin{eqnarray}
{n-1 \choose r-1} p_d^{r-1} &=& \frac{(n-1)!}{(r-1)!(n-r)!} \frac{C_d^{r-1}}{n^{r-1}} \geq \frac{C_d^{r-1}}{(r-1)!} (1-\epsilon)^{r-1}.\label{term1_est}
\end{eqnarray}
The final inequality in~(\ref{term1_est}) is true since
\[\frac{(n-1)!}{(n-r)!n^{r-1}} = \frac{(n-1)\ldots(n-r+1)}{n^{r-1}} \geq \frac{(n-r)^{r-1}}{n^{r-1}} \geq (1-\epsilon)^{r-1}.\] In the above,
the first inequality follows since \(n-k \geq n-r\) for \(1 \leq k \leq r\) and the second inequality follows from the fact that \(r \leq \epsilon n.\)


To estimate the last term in (\ref{t_est4}) we use the fact
\begin{equation}\label{exp_func_est}
1 - x \geq e^{-x - x^2}
\end{equation}
for \(0 < x < \frac{1}{2}.\) For completeness we give a small proof of (\ref{exp_func_est}). For \(0 < x < \frac{1}{2}\) we have
\begin{equation}\label{log_est}
\log(1-x) = - x - R(x)
\end{equation}
where  \(0 \leq R(x) = \sum_{k \geq 2} \frac{x^k}{k}.\) Expanding the series, we have \[R(x) = \frac{x^2}{2} \left(1 + x\frac{2}{3} + x^2\frac{2}{4} + \ldots\right) \leq \frac{x^2}{2}(1+x+x^2 +\ldots) = \frac{x^2}{2(1-x)} \leq x^2\] where the final estimate follows since \(1-x \geq \frac{1}{2}.\) Substituting into (\ref{log_est}) gives~(\ref{exp_func_est}).

We fix \(n_0\) large so that
\begin{equation}\label{p_u_conv}
p_u = \frac{C_u}{n} \leq \frac{C+\alpha_n}{n} \leq \frac{C+1}{n} \leq \frac{1}{2}
\end{equation}
for all \(n \geq n_0.\) The second inequality is true since \(\alpha_n \longrightarrow 0\) as \(n\rightarrow \infty.\) Using the bound (\ref{exp_func_est}) for the last term in (\ref{t_est4}) then gives
\begin{eqnarray}
(1-p_u)^{r(n-r) + {r \choose 2} - r+1} = \left(1-\frac{C_u}{n}\right)^{nr - \frac{r^2}{2} - \frac{3r}{2} + 1} \geq e^{-I_1 -I_2} \label{term2_est1}
\end{eqnarray}
where
\begin{eqnarray} \label{i_1_est}
I_1 = \frac{C_u}{n}\left(nr - \frac{r^2}{2} -\frac{3r}{2} +1\right) \leq \frac{C_u}{n}\left(nr+1\right)  = C_ur + \frac{C_u}{n}
\end{eqnarray}
and
\begin{equation}\label{i_2_est}
I_2 =  \frac{C_u^2}{n^2}\left(nr - \frac{r^2}{2} - \frac{3r}{2} +1\right) \leq \frac{C_u^2}{n^2}\left(nr+1\right) =\frac{C_u^2r}{n^2} + \frac{C_u^2}{n^2} \leq \frac{C_u^2 \epsilon}{n} + \frac{C_u^2}{n^2}.
\end{equation}
The final estimate is true since \(r \leq \epsilon n.\)

Substituting the above two estimates into (\ref{term2_est1}) we have
\begin{eqnarray}
(1-p_u)^{r(n-r) + {r \choose 2} - r+1} &\geq& e^{-C_ur}\exp\left(-\frac{C_u}{2n} - \frac{C_u^2\epsilon}{n}-\frac{C_u^2}{2n^2}\right) \nonumber\\
&\geq& e^{-C_ur}e^{-2\epsilon}\label{term2_est2}
\end{eqnarray}
for all \(n \geq n_1.\) Here \(n_1\) does not depend on the choice of \(r.\) The final estimate holds since we have from (\ref{p_u_conv}) that \(\frac{C_u}{2n} + \frac{C_u^2\epsilon}{n}+\frac{C_u^2}{2n^2}\) converges to zero as \(n \rightarrow \infty.\)


Using (\ref{term2_est2}) and (\ref{term1_est}) in (\ref{t_est4}) gives
\begin{eqnarray}\label{c1_temp_est_aw}
\mathbb{P}(\#{\cal E}_i = r) \geq T_r\frac{C_d^{r-1}}{(r-1)!} (1-\epsilon)^{r-1}e^{-C_ur}e^{-2\epsilon}.
\end{eqnarray}
Using \(C_u = C+\alpha_n\) and \(C_d = C-\alpha_n,\) we have that
\begin{eqnarray}
C_d^{r-1} e^{-C_u r} &=& C^{r-1}(1-\alpha_nC^{-1})^{r-1} e^{-Cr} e^{-\alpha_n r} \nonumber\\
&\geq& C^{r-1}e^{-\alpha_n C^{-1}(r-1) - \alpha_n^2C^{-2} (r-1)} e^{-Cr} e^{-\alpha_n r} \nonumber\\
&=& C^{r-1}e^{-Cr} e^{-\alpha_nC^{-1} (r-1) - \alpha_n^2C^{-2} (r-1)-\alpha_n r} \label{c1_temp_est_wq}
\end{eqnarray}
for all \(n \geq n_2.\) Here \(n_2\) does not depend on the choice of \(r.\) The middle inequality is obtained using (\ref{exp_func_est}) since \(\alpha_n \longrightarrow 0\) as \(n \rightarrow \infty\) and so \(\alpha_n < \frac{1}{2}\) for all \(n\) large.

Fixing \(\omega > 0\) as in  the statement of the Lemma, we also have that \[\max(\alpha_nC^{-1},\alpha_n^2C^{-2},\alpha_n) < \frac{\omega}{3}\] for all \(n \geq N_3.\) Here~\(N_3\) does not depend on the choice of~\(r.\) Thus the exponent in the final term of (\ref{c1_temp_est_wq}) can be bounded as
\begin{equation}
\alpha_nC^{-1}(r-1) + \alpha^2_nC^{-2}(r-1) + \alpha_n r \leq \frac{\omega}{3}(r-1) + \frac{\omega}{3}(r-1) + \frac{\omega}{3}r = \omega r - \frac{2\omega}{3}. \nonumber
\end{equation}
Substituting the above into (\ref{c1_temp_est_wq}) we have
\begin{eqnarray}
C_d^{r-1} e^{-C_u r} \geq C^{r-1}e^{-Cr} e^{-\omega r} e^{\frac{2\omega}{3}}
\end{eqnarray}
and using this in (\ref{c1_temp_est_aw}) gives (\ref{c1_est2}). \(\qed\)



\subsection*{Cross term estimates}
We also need the following estimate on the cross terms to determine the variance in the size of the giant component.
\begin{Lemma}\label{cross}Fix \(C \neq 1\) and fix \(\epsilon,\omega > 0\) such that \(\delta_0\) and \(\delta_1\) defined in~(\ref{delta0_def}) are positive. There is an integer \(N = N(C,\epsilon,\omega)\) so that the following two statements hold. For all \(n \geq N, 1\leq i \neq j \leq n\) and \(1 \leq r_1, r_2 \leq \epsilon n,\) we have
\begin{eqnarray}
&&\mathbb{P}\left(\{\#{\cal E}_i = r_1\} \cap \{ \#{\cal E}_j = r_2\} \cap \{{\cal E}_i \neq {\cal E}_j\}\right) \nonumber\\
&&\;\;\;\;\leq \frac{1}{C}\frac{T_{r_1}e^{-r_1}}{(r_1-1)!}e^{-\delta_1 r_1} \frac{1}{C}\frac{T_{r_2}e^{-r_2}}{(r_2-1)!}e^{-\delta_1 r_2}.\label{cr1_est1}
\end{eqnarray}
Similarly, for all \(n \geq N, 1 \leq i \neq j \leq n\) and \(2 \leq r_1 \leq \epsilon n,\) we have
\begin{eqnarray}
\mathbb{P}\left(\{\#{\cal E}_i = r_1\} \cap \{{\cal E}_i = {\cal E}_j\}\right) \leq 2\epsilon \frac{1}{C}\frac{T_{r_1}e^{-r_1}}{(r_1-1)!}e^{-\delta_0 r_1}.\label{cr1_est2}
\end{eqnarray}
\end{Lemma}
\emph{Proof of Lemma~\ref{cross}}: The proof is analogous to the proof of the upper bound~(\ref{c1_est1}) of Lemma~\ref{c1_est1_lem}. As in the proof of Lemma~\ref{c1_est1_lem}, let \(C_u = C+\alpha_n, C_d = C-\alpha_n\) and \(p_u = \frac{C_u}{n}\) and \(p_d = \frac{C_d}{n}.\)

Also we let \(i = 1\) and \(j = 2\) throughout the proof and the argument holds for all~\(i\) and~\(j.\) Let
\begin{equation}\label{e_def}
E := \{\#{\cal E}_1 = r_1\}\cap \{\#{\cal E}_2 = r_2\} \cap \{{\cal E}_1 \neq {\cal E}_2\}
\end{equation}
and suppose that \({\cal E}_1 = \{1,2,\ldots, r_1\}\) and \({\cal E}_2= \{r_1+1,\ldots,r_1+r_2\}\) for \(r_1,r_2 \geq 2.\)
Let \({\cal T}_{r_1}\) be the set of all trees with vertex set \(\{1,2,\ldots,r_1\}\) and let~\({\cal T}_{r_1,r_2}\) be the set of all trees with vertex set~\(\{r_1+1,\ldots,r_1+r_2\}.\) Since \({\cal E}_1\) has \(\{1,2,\ldots,r_1\}\) as its vertices, there is a random tree \({\cal J}_1 \in {\cal T}_{r_1}\) containing~\(r_1-1\) edges, each of which is open. Moreover, every edge with one endvertex in \({\cal J}_1\) and one endvertex outside, is closed. The number of such edges is \(r_1(n-r_1).\) Similarly, there is a random tree \({\cal J}_2 \in {\cal T}_{r_1,r_2}\) containing \(r_2-1\) edges, each of which is open. Again every edge with one endvertex in \({\cal J}_2\) and one endvertex outside, is closed. In this case, we only need to count the edges not having an endvertex in the tree~\({\cal J}_1.\) The number of such edges is \(r_2(n-r_1-r_2).\)

Let \[E_{r_1,r_2} = E \cap \{{\cal E}_1 = \{1,2,\ldots,r_1\}\} \cap \{{\cal E}_2 = \{r_1+1,\ldots,r_1+r_2\}\}\] and fix trees \({\cal T}_1 \in {\cal T}_{r_1}\) and \({\cal T}_2 \in {\cal T}_{r_1,r_2}.\) We have from the discussion in the previous paragraph that
\begin{equation}
\mathbb{P}\left(E_{r_1,r_2} \cap \{{\cal T}_1 \subseteq {\cal E}_1\} \cap \{{\cal T}_2 \subseteq {\cal E}_2\} \right) \leq A_1 \label{eq_c_bg211}
\end{equation}
where
\begin{equation}\label{a1_def}
A_1 = A_1(r_1,r_2) = p_u^{r_1-1}(1-p_d)^{r_1(n-r_1)}p_u^{r_2-1}(1-p_d)^{r_2(n-r_1-r_2)}.
\end{equation}
In obtaining the above expression, we use the fact that every edge in the complete graph~\(K_n\) is open with probability at most \(p_u\) and closed with probability at most \(1-p_d\) (see (\ref{edge_prob})).
Summing over all possible choices of~\({\cal T}_1\) and~\({\cal T}_2\) we have
\begin{equation}
\mathbb{P}\left(E_{r_1,r_2}\right) \leq T_{r_1}T_{r_2} A_1 \label{eq_c_bg212}
\end{equation}
where as before \(T_r\) denotes the number of labelled trees on \(r\) vertices.

The expression (\ref{eq_c_bg212}) obtains estimates for a particular choice of vertex sets for the components \({\cal E}_1\) and \({\cal E}_2.\) To count the total number of choices for the vertex sets of the components \({\cal E}_1\) and \({\cal E}_2,\) we argue as follows. Since~\(1~\in~{\cal E}_1\) and~\(2 \in {\cal E}_2,\) we choose \(r_1-1\) vertices for~\({\cal E}_1\) and~\(r_2 - 1\) vertices for~\({\cal E}_2\) out of the remaining~\(n-2\) vertices. The number of such choices is
\begin{equation}\label{n_choose}
{n-2 \choose r_1-1,r_2-1} := \frac{(n-2)!}{(r_1-1)!(r_2-1)!(n-r_1-r_2)!}.
\end{equation}
Using (\ref{eq_c_bg212}) and the definition of the event \(E\) in (\ref{e_def}) we therefore have that
\begin{eqnarray}
\mathbb{P}\left(E\right) \leq \frac{(n-2)!}{(r_1-1)!(r_2-1)!(n-r_1-r_2)!}T_{r_1}T_{r_2} A_1. \label{crs_est1}
\end{eqnarray}
Setting \(T_1 = 1,\) we have that the above expression also holds if \(r_1=1\) or \(r_2 = 1.\)
Substituting the expression for \(A_1\) from (\ref{a1_def}) and rearranging terms, we have
\begin{eqnarray}
\mathbb{P}\left(E\right) &\leq& \frac{(n-2)!}{(r_1-1)!(r_2-1)!(n-r_1-r_2)!}T_{r_1}T_{r_2} p_u^{r_1+r_2-2} (1-p_d)^{A_2} \nonumber\\
&=& \frac{T_{r_1}}{(r_1-1)!}\frac{T_{r_2}}{(r_2-1)!} \frac{(n-2)!}{(n-r_1-r_2)!}p_u^{r_1+r_2-2} (1-p_d)^{A_2} \label{e1_def}
\end{eqnarray}
where
\begin{equation}\label{a_2_def}
A_2 = A_2(r_1,r_2) = r_1(n-r_1) + r_2(n-r_1-r_2).
\end{equation}


We evaluate the last two terms in (\ref{e1_def}) separately.
Using \(p_u = \frac{C_u}{n},\) we have
\begin{eqnarray}
\frac{(n-2)!}{(n-r_1-r_2)!} p_u^{r_1+r_2-2} &=& (n-2)\ldots(n-r_1-r_2+1) \frac{C_u^{r_1+r_2-2}}{n^{r_1+r_2-2}} \nonumber\\
&=&C_u^{r_1+r_2-2}\left(1-\frac{2}{n}\right)\ldots \left(1 - \frac{r_1+r_2-1}{n}\right)\nonumber\\
&\leq& C_u^{r_1+r_2-2}.\label{temp112}
\end{eqnarray}

To evaluate the last term in (\ref{e1_def}), we use (\ref{a_2_def}) to get
\begin{equation}\label{temp11}
(1-p_d)^{A_2} \leq \exp\left(-\frac{C_dA_2}{n}\right) = e^{-C_d(r_1+r_2)}\exp\left(\frac{C_d}{n}(r_1^2 + r_1r_2 + r_2^2)\right).
\end{equation}
To obtain the first inequality we use the estimate \(1-x \leq e^{-x}\) with\\\(x= p_d = \frac{C_d}{n}.\)
Substituting (\ref{temp11}) and (\ref{temp112})
into (\ref{e1_def}) we get
\begin{eqnarray}
\mathbb{P}(E) \leq T_{r_1}\frac{C_u^{r_1-1}e^{-C_dr_1}}{(r_1-1)!}T_{r_2}\frac{C_u^{r_2-1}e^{-C_d r_2}}{(r_2-1)!}e^{A_3} \nonumber
\end{eqnarray}
where \[A_3 = \frac{C_d}{n}(r_1^2 + r_1r_2 + r_2^2)  \leq C_d\epsilon(r_1 + r_1 + r_2) \leq 2C_d\epsilon(r_1+r_2).\] The middle estimate is obtained
using \(r_1 \leq \epsilon n\) and \(r_2 \leq \epsilon n.\) Thus
\begin{eqnarray}
\mathbb{P}(E) \leq \frac{T_{r_1}}{(r_1-1)!} \frac{T_{r_2}}{(r_2-1)!}A_4(r_1) A_4(r_2) \label{e_penul}
\end{eqnarray}
where \(A_4(r) = C_u^{r-1}e^{-C_d r}e^{2C_d \epsilon r}.\)

Arguing as in the derivation of (\ref{cross_later}), we have that
\begin{equation}
A_4(r) = C_u^{r-1} e^{-C_d r}e^{2C_d\epsilon r} \leq C^{r-1}e^{-Cr} e^{2C\epsilon r} e^{\omega r}  = \frac{e^{-r}}{C} e^{-\delta_1 r} \label{cross_est_ned}
\end{equation}
for all \(n \geq n_1,\) where \(n_1\) does not depend on the choice of \(r.\) Here \(\omega >0\) is as in the statement of the Lemma and \(\delta_1 > 0\) is as in (\ref{delta0_def}).
Substituting~(\ref{cross_est_ned}) into (\ref{e_penul}) gives (\ref{cr1_est1}).

The proof of (\ref{cr1_est2}) is analogous.
First we write
\[\mathbb{P}\left(\{\#{\cal E}_1 = r_1\} \cap \{{\cal E}_1 = {\cal E}_2\}\right) = \mathbb{P}\left(\{\#{\cal E}_1 = r_1\} \cap \{2 \in {\cal E}_1\}\right).\]
If \(\{\#{\cal E}_1 = r_1\} \cap \{2 \in {\cal E}_1\}\) occurs, then we only need to choose \(r_1-2\) remaining vertices out of the possible \(n-2\) vertices.
The number of such ways is \({n-2 \choose r_1-2}.\) As before, there are \(T_{r_1}\) possible labelled trees with \(r_1\) vertices. Arguing as in the proof of (\ref{c1_est1}),
we have
\begin{eqnarray}
&&\mathbb{P}(\{\#{\cal E}_1 = r_1\} \cap \{2 \in {\cal E}_1\}) \nonumber\\
&&\;\;\;\;\;\;\leq {n-2 \choose r_1-2} T_{r_1} p_u^{r_1-1}(1-p_d)^{r_1(n-r_1)} \nonumber\\
&&\;\;\;\;\;\;=\frac{r_1-1}{n-1} \left({n-1 \choose r_1-1} T_{r_1} p_u^{r_1-1}(1-p_d)^{r_1(n-r_1)} \right). \label{crs_est22}
\end{eqnarray}
The term within the brackets is exactly the term in the right hand side of~(\ref{eq_c_bg2}) with~\(r\) replaced by \(r_1\) and whose estimate is obtained as the right hand side of (\ref{c1_est1}). Also since \(r_1 \leq \epsilon n,\) we have that \(\frac{r_1-1}{n-1} \leq \frac{\epsilon n-1}{n-1} \leq 2\epsilon\) for all~\(n\) large.

Combining the above, we have
\begin{eqnarray} \nonumber
\mathbb{P}(\{\#{\cal E}_1 = r_1\} \cap \{2 \in {\cal E}_1\}) \leq 2\epsilon \frac{1}{C}\frac{T_{r_1}e^{-r_1}}{(r_1-1)!}e^{-\delta_0 r_1}.
\end{eqnarray}
This proves (\ref{cr1_est2}). \(\qed\)

\setcounter{equation}{0}
\renewcommand\theequation{\thesection.\arabic{equation}}
\section{Mid size components estimate}\label{mid_comp}
In Lemmas~\ref{c1_est1_lem} and~\ref{cross} of the previous section, we have obtained probability estimates for non giant components; i.e., components
of size at most \(\epsilon n\) for some constant \(\epsilon > 0.\) The next step is to obtain probability estimates for mid size components; i.e., components whose size lies in the range \([M\log{n}+1,\epsilon n]\) for some constant \(M >0.\) To do so, we use an auxiliary result. Let \(T_1 = 1\) and for integer \(r \geq 2,\) let~\(T_r\) be the number of labelled trees on \(r\) vertices. For \(C > 0,\) define \(q(C)\) as in~(\ref{q_0_def}) and \(\delta  = \delta(C) \geq 0\) be as in~(\ref{del_def2}).
Define
\begin{equation}\label{rem_def}
R_n(C) := \sum_{r = n+1}^{\infty} \frac{T_re^{-r}}{(r-1)!}e^{-\delta r} = \sum_{r = n+1}^{\infty} \frac{T_rC^{r-1}e^{-Cr}}{(r-1)!},
\end{equation}
to be the remainder term of \(q(C).\) We have
\begin{equation}\label{rem_est}
\lim_n R_n(C) = 0.
\end{equation}
We directly use the properties of the random graphs to prove the result.\\
\emph{Proof of~(\ref{rem_est})}:  For \(1 \leq i \leq n,\) we recall that \({\cal E}_i\) denotes the component of the random graph \(G\) containing the vertex \(i\) and so the term
\begin{equation}\label{zn_def}
Z_n(\epsilon) := \sum_{i=1}^{n} \ind(1 \leq \#{\cal E}_i \leq \epsilon n)
\end{equation}
denotes the sum of sizes of all components each of whose size is at most~\(\epsilon n.\)

For \(C \neq 1,\) we fix \(\epsilon,\omega > 0\) so that \(\delta_2  = \delta_2(C,\epsilon,\omega)\) defined in~(\ref{delta0_def}) is positive.
We then use the lower bound (\ref{c1_est2}) in Lemma~\ref{c1_est1_lem} to get
\begin{eqnarray}
\mathbb{E}Z_n(\epsilon) = \sum_{i=1}^{n} \mathbb{P}(1 \leq \#{\cal E}_i \leq \epsilon n) \geq n\sum_{r=1}^{\epsilon n} \frac{1}{C(1-\epsilon)}\frac{T_re^{-r}}{(r-1)!}e^{-\delta_2 r} e^{-2\epsilon + \frac{2\omega}{3}}. \label{use_lw_bd}
\end{eqnarray}
Fix integer \(N \geq 1.\) Choose \(n\) large so that \(\epsilon n \geq N.\) We then have from (\ref{use_lw_bd}) that
\begin{eqnarray}
\sum_{r=1}^{N} \frac{1}{C(1-\epsilon)}\frac{T_re^{-r}}{(r-1)!}e^{-\delta_2 r} e^{-2\epsilon + \frac{2\omega}{3}} \leq  \frac{1}{n}\mathbb{E}Z_n(\epsilon) \leq 1.\label{use_lw_bd2}
\end{eqnarray}
The final estimate is true since \(Z_n(\epsilon) \leq n.\) Also \(\delta_2(C,\epsilon,\epsilon) = \delta -(2C+1)\epsilon\) is a decreasing function of~\(\epsilon\) and is positive for all \(\epsilon> 0\) small. Therefore allowing \(\epsilon = \omega \downarrow 0\)  in~(\ref{use_lw_bd2}) gives
\begin{eqnarray}
\sum_{r=1}^{N} \frac{1}{C}\frac{T_re^{-r}}{(r-1)!}e^{-\delta r}  \leq  1\label{use_lw_bd3}
\end{eqnarray}
where \(\delta > 0\) is as in (\ref{del_def2}).
Allowing \(N \rightarrow \infty\) gives \(q(C) \leq 1\) for all \(C \neq 1.\)

For \(C = 1,\) we again use (\ref{use_lw_bd3}). Allowing \(C \uparrow 1\) in (\ref{use_lw_bd3}) gives
\begin{eqnarray}
\sum_{r=1}^{N} \frac{T_re^{-r}}{(r-1)!} \leq  1. \label{use_lw_bd4}
\end{eqnarray}
Again allowing \(N \rightarrow \infty\) gives \(q(1) \leq 1.\) Thus
\begin{equation}\label{qc_less_one}
q(C) \leq 1
\end{equation}
for all \(C > 0\) and this proves~(\ref{rem_est}).\(\qed\)

The following Lemma establishes that with high probability there are no components with size greater than a constant multiple of \(\log{n}\) and less than~\(\epsilon n.\)

For any \(M  > 0 \) and \(\epsilon > 0,\) let
\begin{equation}\label{b_m_def}
B(M,\epsilon) := \bigcup_{1 \leq i \leq n}\{M \log{n}+1 \leq \#{\cal E}_i \leq \epsilon n\}
\end{equation}
denote the event that there is a component whose size lies in the range \([M\log{n}+1,\epsilon n].\) Let \(\delta=  \delta(C)\) be as in (\ref{del_def2}). We have the following Lemma.
\begin{Lemma}\label{lem_mid} Fix \(C \neq 1\)
and \(\epsilon > 0\) so that~\(\delta_0\) defined in (\ref{delta0_def}) is positive. Fix \(M > \frac{1}{\delta_0}\) and \(\theta < M\delta_0 - 1.\) There is a positive constant \(N=N(C, \epsilon,M, \theta)\) such that
\begin{equation}\label{bm_est}
\mathbb{P}(B(M,\epsilon)) \leq \frac{1}{n^{\theta}}
\end{equation}
for all \(n \geq N.\)
\end{Lemma}
\emph{Proof of Lemma~\ref{lem_mid}}: For \(M \log{n}+1 \leq r \leq \epsilon n\) we have that
\[e^{-\delta_0 r} \leq \exp\left(-M \delta_0 \log{n}\right) = \frac{1}{n^{M\delta_0}}.\]
Therefore from (\ref{c1_est1}) we have for any \(1 \leq i \leq n\) that
\begin{equation} \nonumber
\mathbb{P}(\#{\cal E}_i = r) \leq \frac{T_re^{-r}}{C(r-1)!}e^{-\delta_0 r}\leq \frac{T_re^{-r}}{C(r-1)!}\frac{1}{n^{M\delta_0}}.
\end{equation}
Hence
\begin{eqnarray}
\mathbb{P}(M \log{n}+1 \leq \#{\cal E}_i \leq \epsilon n) &\leq& \frac{1}{n^{M\delta_0}}\sum_{r = M\log{n}+1}^{\epsilon n} \frac{T_re^{-r}}{C(r-1)!} \nonumber\\
&\leq& \frac{1}{Cn^{M\delta_0}}R_{M\log{n}}(1) \nonumber
\end{eqnarray}
where \(R_{M\log{n}}(.)\) is the remainder term defined in (\ref{rem_def}).
Using (\ref{rem_est}) we have that \(R_{M\log{n}}(1) \leq C\) for all \(n\) large. Thus we have
\begin{equation}\nonumber
\mathbb{P}(M \log{n} +1 \leq \#{\cal E}_i \leq \epsilon n) \leq \frac{1}{n^{M\delta_0}}
\end{equation}
for all \(n \geq n_0.\) Here \(n_0\) does not depend on the choice of \(i.\)

Using the above bound, we have
\begin{eqnarray}
\mathbb{P}(B(M,\epsilon)) \leq \sum_{i=1}^{n}\mathbb{P}(M\log{n} +1 \leq \#{\cal E}_i \leq \epsilon n) \leq \frac{n}{n^{M\delta_0}}. \nonumber
\end{eqnarray}
If \(\theta < M\delta_0 -1,\) then the right hand size is at most \(\frac{1}{n^{\theta}}\) for all \(n\) large.\(\qed\)




\setcounter{equation}{0}
\renewcommand\theequation{\thesection.\arabic{equation}}
\section{Mean and Variance estimates for \(Z_n(\epsilon)\)}\label{mean_zn}
In this section, we obtain mean and variance estimates for the term \(Z_n(\epsilon)\) defined in (\ref{zn_def}), representing the number of vertices present in non giant components, i.e. components whose size is at most \(\epsilon n.\)

We first have some preliminary estimates. Fix \(\epsilon,\omega > 0\) such that \(\delta_i = \delta_i(C,\epsilon,\omega)\) defined in (\ref{delta0_def}) is positive. By definition \(\delta_i(C,\eta_1,\eta_2)\) is positive for all \(0 < \eta_1 < \epsilon\) and \(0 < \eta_2 < \omega.\)
For \(i = 0,1,2,\) let
\begin{eqnarray}
q_i(C,\epsilon,\omega) &:=& e^{-C} + \sum_{r=2}^{\infty} \frac{1}{C}\frac{T_re^{-r}}{(r-1)!}e^{-\delta_i r} \label{q_i_def}
\end{eqnarray}
We have the following result.
Fix \(C \neq 1\) and \(i \in \{0,1,2\}.\) We have
\begin{equation}\label{q_i_lim}
\lim_{\epsilon \downarrow 0} q_i(C,\epsilon,\epsilon) = q(C).
\end{equation}
\emph{Proof of~(\ref{q_i_lim})}: We first prove for \(i = 2.\)
Let \(\epsilon_m \downarrow 0\) be any sequence.
Since \((1-\epsilon)e^{\epsilon} \leq 1\)  for any \(\epsilon > 0,\) we have
\begin{eqnarray}\nonumber
q_2(C,\epsilon_m,\epsilon_m) = \sum_{r=1}^{\infty} b_r(1-\epsilon_m)^{r}e^{\epsilon_m r} \leq \sum_{r=1}^{\infty}b_r = q(C)
\end{eqnarray}
where \(b_r = \frac{T_rC^{r-1}e^{-Cr}}{(r-1)!}\) for \(r \geq 1.\)
Therefore \(\limsup_m q_2(C,\epsilon_m,\epsilon_m) \leq q(C).\) We obtain the lower bound as follows. For any fixed integer \(N \geq 1\) we have
\begin{eqnarray}\nonumber
q_2(C,\epsilon_m,\epsilon_m) = \sum_{r=1}^{\infty} b_r(1-\epsilon_m)^{r}e^{\epsilon_m r} \geq \sum_{r=1}^{N}b_r(1-\epsilon_m)^{r} e^{\epsilon_m r}
\end{eqnarray}
Allowing \(m \rightarrow \infty\) and then \(N \rightarrow \infty\) gives
\begin{equation}\label{low_bdq2}
\liminf_m q_2(C,\epsilon_m,\epsilon_m) \geq q(C).
\end{equation}
Therefore \(\lim_m q_2(C,\epsilon_m,\epsilon_m) = q(C)\) and (\ref{q_i_lim}) holds for \(i  = 2.\)

We now see that (\ref{q_i_lim}) holds for \(i = 1\) and an analogous proof holds for~\(i = 0.\) We again let \(\epsilon_m \downarrow 0\) be any sequence and
obtain
\begin{eqnarray}\nonumber
q_1(C,\epsilon_m,\epsilon_m) = \sum_{r=1}^{\infty} b_re^{\beta_m r} \geq \sum_{r=1}^{\infty}b_r = q(C)
\end{eqnarray}
where \(\beta_m = 2C\epsilon_m + \epsilon_m.\) Therefore \(\liminf_m q_1(C,\epsilon_m,\epsilon_m) \geq q(C).\)
To obtain the upper bound, we use the fact that
\(\delta_1(C,\epsilon_m,\epsilon_m) = \delta - \beta_m > 0\) for all \(m \geq M_0\) large. Fix \(m \geq M_0\) and any integer \(N \geq 1.\) We then have that
\begin{eqnarray}\label{q1_up_bd}
q_1(C,\epsilon_m,\epsilon_m) = \sum_{r=1}^{\infty} b_re^{\beta_m r} = \sum_{r=1}^{N}b_re^{\beta_m r} + R_{1,N}(C,\epsilon_m)
\end{eqnarray}
where
\begin{equation}\label{rem111}
R_{1,N}(C,\epsilon_m) := \sum_{r = N+1}^{\infty}b_re^{\beta_m r} = \sum_{r=N+1}^{\infty} \frac{T_re^{-r}}{(r-1)!}e^{-(\delta-\beta_m) r} \leq \sum_{r=N+1}^{\infty} \frac{T_re^{-r}}{(r-1)!}
\end{equation}
and the last term is precisely \(R_N(1)\) defined in (\ref{rem_def}).

Substituting the above into (\ref{q1_up_bd}) we have
\begin{eqnarray}\nonumber
q_1(C,\epsilon_m,\epsilon_m) \leq \sum_{r=1}^{N}b_re^{\beta_m r} + R_N(1)
\end{eqnarray}
for all \(m \geq M_0\) and for any fixed integer \(N \geq 1.\)
Allowing \(m \rightarrow \infty\) in the above and using \(\beta_m \longrightarrow 0,\) we have
\[\limsup_m q_1(C,\epsilon_m,\epsilon_m) \leq \sum_{r=1}^{N}b_r +  R_N(1).\] Allowing \(N \rightarrow \infty\) in the above and using (\ref{rem_est}), we get
\[\limsup_m q_1(C,\epsilon_m,\epsilon_m) \leq \sum_{r \geq 1} b_r = q(C).\] Thus \(\lim_m q_1(C,\epsilon_m,\epsilon_m) = q(C)\) and this proves (\ref{q_i_lim}) for \(i = 1.\)~\(\qed\)

\subsubsection*{Mean estimates for \(Z_n(\epsilon)\)}
We have the following bounds on the mean of
\begin{equation}\label{x_def}
Z_n(\epsilon) = \sum_{1 \leq i \leq n} \ind(1 \leq {\cal E}_i \leq \epsilon n).
\end{equation}
For \(C > 0,\) let \(q = q(C)\) be
as defined in~(\ref{q_0_def}).
\begin{Lemma}\label{exp_x_lem} Fix \(0 < \gamma < 1\) and \( C \neq 1.\) There is a positive constant \(\epsilon_1 = \epsilon_1(C,\gamma)\) so that the following statement holds for all \(0 < \epsilon < \epsilon_1.\) There is a positive constant \(N = N(C,\epsilon, \gamma)\) so that for all \(n \geq N,\) we have
\begin{equation}\label{exp_x_est}
nq(C)(1-\gamma) \leq \mathbb{E}Z_n(\epsilon) \leq nq(C)(1+\gamma).
\end{equation}
\end{Lemma}
\emph{Proof of Lemma~\ref{exp_x_lem}}: Fix \(i \in \{ 0,2\}\) and let \(q_i(.,.,.)\) be as defined in (\ref{q_i_def}). From (\ref{q_i_lim}), we have that \(q_i(C,\eta,\eta) \longrightarrow q(C)\) as \(\eta \rightarrow 0.\) For any fixed \(\eta > 0,\) the term \(\delta_2 = \delta_2(C,\eta,\eta) = \delta -\log(1-\eta) + \eta\) is positive. Here \(\delta = \delta(C) > 0\) is as in (\ref{del_def2}). Also, the term \(\delta_0 = \delta_0(C,\eta,\eta) =\delta - (2C+1)\eta >0\) for all \(0 < \eta < \frac{\delta}{2C+1}.\)

For a fixed \(0 < \gamma < 1,\) we therefore let \(\epsilon_1  = \epsilon_1(C,\gamma) > 0\) be small so that for all \(0 < \epsilon < \epsilon_1,\) the following statements hold.
The term
\begin{equation}\label{eq_eps_wq}
\left(1-\frac{3\gamma}{4}\right)\frac{e^{-\frac{4\epsilon}{3}}}{(1-\epsilon)}  \geq 1-\gamma.
\end{equation}
The terms \(\delta_1(C,\epsilon,\epsilon)\) and \(\delta_2(C,\epsilon,\epsilon)\) are positive and
\begin{equation}\label{q_bd02}
q_0(C,\epsilon,\epsilon)\leq q(C) (1+\gamma)  \text{ and } q_2(C,\epsilon,\epsilon)\geq q(C)\left(1 - \frac{\gamma}{2}\right).
\end{equation}

Fix an \(\epsilon > 0\) so that the above statements hold. We obtain an upper bound for \(Z_n(\epsilon)\) first.
Using the upper bound (\ref{c1_est1}) in Lemma~\ref{c1_est1_lem}, we have that
\begin{eqnarray}
\frac{1}{n}\mathbb{E}Z_n(\epsilon) = \frac{1}{n}\sum_{i=1}^{n}\mathbb{P}(1 \leq \#{\cal E}_i \leq \epsilon n) \leq \sum_{r=1}^{\epsilon n} \frac{1}{C}\frac{T_re^{-r}}{(r-1)!}e^{-\delta_0 r} \leq q_0(C,\epsilon,\epsilon).\;\;\;\label{exp_x_est2}
\end{eqnarray}
Using (\ref{q_bd02}), we have
\begin{equation}\label{exp_x_up}
\frac{1}{n} \mathbb{E}Z_n(\epsilon) \leq q(C)(1+\gamma).
\end{equation}

To obtain a lower bound on \(\mathbb{E}Z_n(\epsilon),\) we use the lower bound (\ref{c1_est2}) in Lemma~\ref{c1_est1_lem} with \(\epsilon = \omega\) to get
\begin{eqnarray}
\frac{1}{n}\mathbb{E}Z_n(\epsilon) &=& \frac{1}{n}\sum_{i=1}^{n}\mathbb{P}(1 \leq \#{\cal E}_i \leq \epsilon n) \nonumber\\
&\geq& \sum_{r=1}^{\epsilon n} \frac{1}{C(1-\epsilon)}\frac{T_re^{-r}}{(r-1)!}e^{-\delta_2 r} e^{-\frac{4\epsilon}{3}} \nonumber\\
&=& \left(q_2(C,\epsilon,\epsilon)- R_{2,n}(C,\epsilon)\right)\frac{e^{-\frac{4\epsilon}{3}}}{(1-\epsilon)}\label{exp_x_est11}
\end{eqnarray}
where
\begin{equation}\nonumber
0 \leq R_{2,n}(C,\epsilon) = \sum_{r=\epsilon n+1}^{\infty} \frac{T_rC^{r-1}e^{-Cr}}{(r-1)!}e^{-\delta_2 r} \leq \sum_{r=\epsilon n+1}^{\infty} \frac{T_rC^{r-1}e^{-Cr}}{(r-1)!} = R_{\epsilon n}(C).
\end{equation}
Here \(R_n(.)\) is the remainder term as in (\ref{rem_def}).

Since the remainder \(R_{\epsilon n}(C) \longrightarrow 0\) as \(n \rightarrow \infty\) (see (\ref{rem_est})), we have
\[R_{\epsilon n}(C) \leq q(C)\frac{\gamma}{4}\] for all large \(n.\) Substituting the above and the estimate for \(q_2(C,\epsilon,\epsilon)\)~(see~(\ref{q_bd02})) into (\ref{exp_x_est11}) we have
\begin{equation}\label{exp_x_down}
\frac{1}{n} \mathbb{E}Z_n(\epsilon) \geq q(C)\left(1-\frac{3\gamma}{4}\right)\frac{e^{-\frac{4\epsilon}{3}}}{(1-\epsilon)} \geq q(C)(1-\gamma)
\end{equation}
by our choice of \(\epsilon > 0\) in (\ref{eq_eps_wq}).
Combining (\ref{exp_x_up}) and (\ref{exp_x_down}) gives (\ref{exp_x_est}).\(\qed\)


\subsubsection*{Variance estimate for \(Z_n(\epsilon)\)}
We have the following estimate on the variance of \(Z_n(\epsilon).\)
\begin{Lemma}\label{var_x_lem} Fix \(0 < \gamma < 1\) and \( C \neq 1.\) There is a positive constant \(\epsilon_2 = \epsilon_2(C,\gamma)\) so that the following statement holds for all \(0 < \epsilon < \epsilon_2.\) There is a positive integer \(N = N(C,\epsilon, \gamma)\) so that for all \(n \geq N,\) we have
\begin{eqnarray}\label{var_x_est}
var(Z_n(\epsilon)) \leq nq(C)(1+\gamma) + 4n^2q^2(C)\gamma \nonumber
\end{eqnarray}
where \(q(C)\) is as in (\ref{q_0_def}).
\end{Lemma}
\emph{Proof of Lemma~\ref{var_x_lem}}: For \(i = 0,1,2\) let \(\delta_i = \delta_i(C,\epsilon,\epsilon).\)
For \(1 \leq i \leq n,\) let \(X_i = \ind(1 \leq \#{\cal E}_i \leq \epsilon n)\) denote the indicator function of the event that the component containing vertex \(i\) has size at most \(\epsilon n.\) From (\ref{zn_def}) we have that \(Z_n(\epsilon) = \sum_{i=1}^{n} X_i\) and so
\begin{equation}
\mathbb{E}Z^2_n(\epsilon) = \sum_{i=1}^{n}\sum_{j=1}^{n}\mathbb{E}X_iX_j = \sum_{i=1}^{n}\mathbb{E}X^2_i + \sum_{i \neq j} \mathbb{E}(X_iX_j) . \label{var_est1}
\end{equation}
Since each \(X_i\) is an indicator function, we have that \(X_i^2 = X_i\) and so the first term in (\ref{var_est1}) is \(\sum_{i=1}^{n} \mathbb{E}X_i = \mathbb{E}Z_n(\epsilon).\) We evaluate the second term in (\ref{var_est1}) using the cross term estimates in Lemma~\ref{cross}.
For \(1 \leq i \neq j \leq n\) we write
\begin{equation}\label{crs_term1}
E(X_iX_j) =  \mathbb{E}(X_iX_j\ind({\cal E}_i = {\cal E}_j) + \mathbb{E}(X_iX_j\ind({\cal E}_i \neq {\cal E}_j))
\end{equation}
and evaluate each term separately.

Using the cross term estimate~(\ref{cr1_est1}) of Lemma~\ref{cross} we have
\begin{eqnarray}
\mathbb{E}(X_iX_j\ind({\cal E}_i \neq {\cal E}_j)) &=& \sum_{r_1=1}^{\epsilon n}\sum_{r_2=1}^{\epsilon n} \mathbb{P}\left(\{\#{\cal E}_i = r_1\}\cap \{\#{\cal E}_j = r_2\} \cap \{{\cal E}_i \neq {\cal E}_j\}\right)\nonumber\\ &\leq& \sum_{r_1=1}^{\epsilon n}\sum_{r_2=1}^{\epsilon n} \frac{1}{C}\frac{T_{r_1}e^{-r_1}}{(r_1-1)!}e^{-\delta_1 r_1} \frac{1}{C}\frac{T_{r_2}e^{-r_2}}{(r_2-1)!}e^{-\delta_1 r_2} \nonumber\\
&\leq& q^2_1(C,\epsilon,\epsilon) \label{tempax1}
\end{eqnarray}
where \(q_1(.,.,.)\) is as in (\ref{q_i_def}). 

Using (\ref{cr1_est2}) of Lemma~\ref{cross}, we similarly have
\begin{eqnarray}
\mathbb{E}(X_iX_j\ind({\cal E}_i = {\cal E}_j)) &=& \sum_{r_1=2}^{\epsilon n}\mathbb{P}\left(\{\#{\cal E}_i = r_1\} \cap \{{\cal E}_i = {\cal E}_j\}\right)\nonumber\\
&\leq& \sum_{r_1=2}^{\epsilon n}2\epsilon \frac{1}{C}\frac{T_{r_1}e^{-r_1}}{(r_1-1)!}e^{-\delta_0 r_1} \nonumber\\
&\leq& 2\epsilon q_0(C,\epsilon,\epsilon) \label{tempax2}
\end{eqnarray}

Combining (\ref{tempax1}) and (\ref{tempax2}) and substituting in (\ref{crs_term1}) we have
\begin{equation}\label{crs_term2}
E(X_iX_j) \leq q^2_1(C,\epsilon,\epsilon) + 2\epsilon q_0(C,\epsilon,\epsilon).
\end{equation}
The right hand side of (\ref{crs_term2}) converges to \(q^2(C)\) as \(\epsilon \rightarrow 0\) using (\ref{q_i_lim}). Here \(q(C)\) is as defined in (\ref{q_0_def}). Fix \(\epsilon_2 = \epsilon_2(C,\gamma) > 0\) small so that \(\epsilon_2 < \epsilon_1(C,\gamma)\) and
the right hand side of (\ref{crs_term2}) is at most \(q^2(C)(1+\gamma)^2.\) Here \(\epsilon_1(.,.)\) is as defined in the statement of Lemma~\ref{exp_x_lem}. 

Fixing \(0 < \epsilon < \epsilon_2\) we therefore have \(\sum_{i \neq j} \mathbb{E}X_iX_j \leq n(n-1) q^2(C)(1+\gamma).\) This evaluates the second term in (\ref{var_est1}). To evaluate the first term in (\ref{var_est1}), we use the fact that each \(X_i, 1 \leq i \leq n\) is an indicator function. Therefore \(X_i^2 = X_i\) and the first term in (\ref{var_est1}) is \(\sum_{i=1}^{n} \mathbb{E}X_i = \mathbb{E}Z_n(\epsilon).\)
Using the above estimates for the two terms of (\ref{var_est1}) we have
\begin{equation}
\mathbb{E}Z^2_n(\epsilon) \leq \mathbb{E}Z_n(\epsilon) + n(n-1)q^2(C)(1+\gamma)^2 \leq nq(C)(1+\gamma) + n^2q^2(C)(1+\gamma)^2 \nonumber
\end{equation}
where the final estimate is obtained using the upper bound in (\ref{exp_x_est}) of\\Lemma~\ref{exp_x_lem}. This is possible by our choice of \(\epsilon > 0.\) Using the lower bound in (\ref{exp_x_est}), we therefore have
\begin{eqnarray}
var(Z_n(\epsilon))  &=& \mathbb{E}Z^2_n(\epsilon) - (\mathbb{E}Z_n(\epsilon))^2 \nonumber\\
&\leq& nq(C)(1+\gamma) + n^2q^2(C)(1+\gamma)^2 - n^2q^2(C)(1-\gamma)^2 \nonumber\\
&=& nq(C)(1+\gamma) + 4n^2q^2(C)\gamma. \nonumber
\end{eqnarray}
This proves the lemma. \(\qed\)




\setcounter{equation}{0}
\renewcommand\theequation{\thesection.\arabic{equation}}
\section{Subcritical case}\label{pf_sub}

For the subcritical case, there are two proofs. If we assume property \((ii)\) of Theorem~\ref{th_q0} that \(q(C)  = 1\) for \(0 < C < 1,\) we then obtain a weaker version of Theorem~\ref{th2} using Chebychev's inequality and the mean and variance estimates of Lemmas~\ref{exp_x_lem} and~\ref{var_x_lem}. In Section~\ref{pf_sup}, we provide such a proof for the supercritical case~\(C > 1.\)

In what follows we prove Theorem~\ref{th2} \emph{without assuming} Theorem~\ref{th_q0}. We use the result of Theorem~\ref{th2} to give a probabilistic proof of properties \((i)-(iii)\) of Theorem~\ref{th_q0} in the next Section.

Fixing \(C < 1,\) we have from Lemma~\ref{lem_mid} that with high probability there are no mid size components whose size lies in the range \([M\log{n}+1,\epsilon n].\) To obtain decay for giant components, i.e., components of size larger than \(\epsilon n,\) we proceed as follows: We first show that with high probability, all vertices in the component \({\cal E}_i\) containing the vertex~\(i, 1 \leq i \leq n,\) are within a distance of order of \(\log{n}\) from the vertex \(i\) and then see that it is not possible to contain order of~\(n\) vertices within such short distance.

For integer \(t \geq 1,\) let \({\cal N}_t(i)\) denote the set of vertices at a distance of \(t\) from vertex \(i\) in the random graph \(G.\) Therefore \(\#{\cal N}_t(i)\) denotes the number of vertices at a distance~\(t\) from vertex~\(i.\) Here distance between vertices~\(x\) and~\(y\) refers to the graph distance and is the number of edges in the path with the least number of edges between~\(x\) and~\(y.\) Since there are at most~\(n\) vertices in~\({\cal E}_i,\) we have that~\(1 \leq t \leq n.\) Recalling the definition of the sequence~\(\alpha_n\) from~(\ref{edge_prob}), we have
the following Lemma.
\begin{Lemma}\label{lem_diam} Fix \(C < 1\) and let \(C_u = C+\alpha_n.\) There is an integer \(N \geq 1\) so that for all \(n \geq N\) we have that \(C_u < 1.\) For all \(n \geq N, 1 \leq i \leq n\) and \(1 \leq t \leq n,\) we have
\begin{equation}\label{diam_exp}
\mathbb{E}(\#{\cal N}_t(i)) \leq C_u^{t}
\end{equation}
and 
\begin{equation}\label{diam_var}
\mathbb{E}\left((\#{\cal N}_t(i))^2\right) \leq \frac{C_u^{t}}{1-C_u}.
\end{equation}
\end{Lemma}
\emph{Proof of Lemma~\ref{lem_diam}}: We first prove (\ref{diam_exp}). We set \(i = 1\) throughout and for \(1 \leq t \leq n\) define \({\cal N}_t := {\cal N}_t(1).\) The proof holds for arbitrary \(1 \leq i \leq n.\)
Let \(S_0  = \{1\}\) and \(U_0 = \{1,2,\ldots,n\}\setminus S_0.\)
Fixing \(t \geq 1\) we have
\begin{eqnarray}\label{exp_nt1}
\mathbb{E}(\#{\cal N}_t) = \sum_{S_0=\{1\}, S_1,\ldots, S_{t-1}}\mathbb{E}\left(\#{\cal N}_t \ind({\cal N}_{t-1} = S_{t-1},\ldots,{\cal N}_1 = S_1, {\cal N}_0 = S_0)\right)
\end{eqnarray}
where the summation is over all subsets \(S_1,\ldots S_{t-1},\) of \(\{1,2,\ldots,n\}.\)
For fixed \(S_1,\ldots,S_{t-1},\) define the event
\begin{equation}\label{f_t_def}
F_{t-1} := \{{\cal N}_{t-1} = S_{t-1},\ldots,{\cal N}_1 = S_1, {\cal N}_0=S_0\}
\end{equation}
and the set
\begin{equation}\label{u_t_def}
U_{t-1} = \{1,2,\ldots,n\} \setminus\left(\cup_{i=0}^{t-1} S_i\right).
\end{equation}
If the sets \(\{S_i\}\) are not mutually disjoint, we have that \(\ind(F_{t-1}) = 0.\)
If however the indicator function \(\ind(F_{t-1}) = 1,\) then the vertices in~\({\cal N}_t\) satisfy the following properties:\\
\((i)\) If \(v \in {\cal N}_t,\) then \(v\) is adjacent (i.e. connected by an open edge) to some vertex in \(S_{t-1}.\) Also, the vertex \(v\) is not adjacent to any vertex in~\(S_i\) for~\(0 \leq i \leq t-1.\)\\
\((ii)\) All vertices of \({\cal N}_t\) are in \(U_{t-1}.\) \\
\emph{Proof of \((i)-(ii)\)}: The first statement of \((i)\) is true as follows. Fix \(v \in {\cal N}_{t}.\) There is a path \((e(1,i_1),e(i_1,i_2),\ldots,e(i_{t-1},i_t = v))\) consisting of~\(t\) open edges from the vertex~\(1\) to vertex~\(v.\) Here we represent the open edge joining the vertices \(i\) and \(j\) as \(e(i,j).\) The vertex \(i_{t-1}\) is at a distance of~\(t-1\) from vertex~\(1\) and therefore~\(i_{t-1} \in S_{t-1}.\) For the second statement, we argue as follows. If the vertex \(v\) is adjacent to some vertex in \(S_i\) for some \( 1 \leq i \leq t-2,\) then the distance between the vertex \(1\) and the vertex \(v\) is at most~\(t-1.\) This is a contradiction since the vertex \(v\) is at a distance of~\(t\) from vertex~\(1.\)

For proving \((ii),\) is true, we suppose that \(v \notin U_{t-1}.\) From the definition~(\ref{u_t_def}) for \(U_{t-1},\) we then have that \(v \in \cup_{i=0}^{t-1} S_i,\) a contradiction to property~\((i)\) proved above. \(\qed\)



From properties \((i)-(ii)\) above we have that
\begin{eqnarray}\label{nt_est1}
\#{\cal N}_t\ind(F_{t-1}) \leq \sum_{y \in S_{t-1}}\sum_{z \in U_{t-1}}X_{y,z}\ind(F_{t-1}) = \sum_{y \in S_{t-1}} J_y \ind(F_{t-1})
\end{eqnarray}
where \(X_{y,z}\) denotes the indicator function of the event that the edge between vertices \(y\) and \(z\) is open and \(J_y = \sum_{z \in U_{t-1}}X_{y,z}\) for \(y \in S_{t-1}.\) We have an upper bound in (\ref{nt_est1}) since a single vertex \(z \in U_{t-1}\) can be connected to multiple vertices in \(S_{t-1}.\)
Also we have that the following property.
\begin{eqnarray}
&&\text{ The event \(F_{t-1}\) defined in (\ref{f_t_def}) is independent of }\nonumber\\
&&\;\;\;\;\text{ the indicator functions~\(\{X_{y,z}\}_{y \in S_{t-1}, z \in U_{t-1}}.\)}\label{indep_prop}
\end{eqnarray}
\emph{Proof of (\ref{indep_prop})}: For subsets \(A, B \subseteq \{1,2,\ldots,n\},\) let \((A,B)\) denote the set of all edges with one endvertex in \(A\) and other endvertex in \(B.\) For integer~\(i \geq 0,\) define the event \[V_i = \bigcap_{e \in (S_i,S_{i+1})} \{e \text{ is open }\} \bigcap \bigcap_{e \in (S_i,U_{i}\setminus S_{i+1})} \{e \text{ is closed}\}.\] The event \(V_i\) depends only on the set of edges having an endvertex in the set~\(S_i.\)

The event \(F_1 = \{N_1 = S_1, N_0 = S_0\}\) can be written as \(F_1 = V_0\) and so the event~\(F_1\) depends only on the state
of edges containing \(S_0= \{1\}\) as the endvertex. Similarly, the event \[F_2 = \{N_2 = S_2\} \cap F_1 = V_1 \cap V_0\] depends only on the state of the edges that have an endvertex in \(S_0 \cup S_1.\) In particular, the event \(F_2\) does not depend on the state of edges having both endvertices in \(U_1 = \{1,2,\ldots,n\}\setminus (S_0\cup S_1).\) Since \(S_2\subset U_1\) and  \(U_2 \subset U_1,\) any edge with one endvertex in \(S_2\) and other endvertex in \(U_2\) has both endvertices in \(U_1.\) Therefore the event \(F_2\) does not depend on the state of edges having one endvertex in~\(S_2\) and other endvertex in~\(U_2.\) Continuing this way inductively,  the event \(F_{t-1}  = \{N_{t-1} = S_{t-1} \} \cap F_{t-2}\) does not depend on the state of edges having one endvertex in \(S_{t-1}\) and other endvertex in~\(U_{t-1}.\)\(\qed\)

Setting \(p_u = \frac{C_u}{n}\) we have from (\ref{nt_est1}) and (\ref{indep_prop}) that
\begin{equation}
\mathbb{E}\#{\cal N}_t\ind(F_{t-1}) \leq \sum_{y \in S_{t-1}}\mathbb{E}J_y\ind(F_{t-1}) = \sum_{y \in S_{t-1}}\mathbb{E}J_y\mathbb{P}(F_{t-1}) \leq C_u\#S_{t-1} \mathbb{P}(F_{t-1})\label{nt_est2}
\end{equation}
where \(J_y = \sum_{z \in U_{t-1}}X_{y,z}\) is as defined in (\ref{nt_est1}).
The final inequality is obtained because, for any fixed \(y \in S_{t-1}\) we have that
\begin{equation}\label{mean_jy}
\mathbb{E}(J_y) = \sum_{z \in U_{t-1}} \mathbb{E}X_{y,z} \leq \#U_{t-1}p_u \leq np_u = C_u.
\end{equation}

Substituting (\ref{nt_est2}) into (\ref{exp_nt1}) we have
\begin{eqnarray}
\mathbb{E}(\#{\cal N}_t) &\leq& C_u \sum_{S_0 = \{1\},S_1,\ldots, S_{t-1}} \#S_{t-1} \mathbb{P}(F_{t-1}) \nonumber\\
&=& C_u \sum_{S_0 = \{1\},S_1,\ldots, S_{t-1}} \#S_{t-1} \mathbb{P}({\cal N}_{t-1} = S_{t-1},\ldots {\cal N}_1 = S_1, {\cal N}_0=S_0) \nonumber\\
&=& C_u \mathbb{E}(\#{\cal N}_{t-1}). \label{rec_eq}
\end{eqnarray}
Continuing this iteratively we get \(\mathbb{E}(\#{\cal N}_t) \leq C_u^{t}\mathbb{E}(\#{\cal N}_0) = C_u^t.\) This proves~(\ref{diam_exp}).

The proof of (\ref{diam_var}) is analogous. Indeed proceeding as before, we have
\begin{equation}
\mathbb{E}(\#{\cal N}_t)^2 = \sum_{S_0=\{1\},S_1,\ldots, S_{t-1}}\mathbb{E}\left( \left(\#{\cal N}_t\right)^2 \ind({\cal N}_{t-1} = S_{t-1},\ldots,{\cal N}_1 = S_1, {\cal N}_0=S_0)\right).\\
\label{exp_nt2}
\end{equation}
Defining \(F_{t-1}\) as in (\ref{f_t_def}) and using (\ref{indep_prop}), we have
\begin{eqnarray}
\mathbb{E}(\#{\cal N}_t)^2\ind(F_{t-1}) \leq \mathbb{E}\left(\sum_{y \in S_{t-1}}J_y\right)^2 \ind(F_{t-1})= \mathbb{E}\left(\sum_{y \in S_{t-1}}J_y\right)^2\mathbb{P}(F_{t-1}) \label{tempax_124}
\end{eqnarray}
where \(J_y = \sum_{z \in U_{t-1}} X_{y,z}\) is as defined in (\ref{nt_est1}).

If \(y_1 \neq y_2,\) the random variables \(J_{y_1}\) and \(J_{y_2}\) are independent
and so
\begin{eqnarray}
\mathbb{E}\left(\sum_{y \in S_{t-1}}J_y\right)^2 &=& var\left(\sum_{y \in S_{t-1}}J_y\right) + \left(\mathbb{E}\sum_{y \in S_{t-1}}J_y\right)^2 \nonumber\\
&=& \sum_{y \in S_{t-1}}var(J_y)  + \left(\sum_{y \in S_{t-1}} \mathbb{E}J_y\right)^2. \label{tempax_123}
\end{eqnarray}
The mean \(\mathbb{E}J_y \leq C_u\) using (\ref{mean_jy}) and regarding the variance, we have that
\begin{equation}\label{var_y}
var(J_y)  = \sum_{z\in U_{t-1}}var(X_{y,z}) \leq \sum_{z \in U_{t-1}} \mathbb{E}X^2_{y,z} = \sum_{z \in U_{t-1}} \mathbb{E}X_{y,z} \leq C_u.
\end{equation}
The first equality in (\ref{var_y}) holds since the random variables \(\{X_{y,z}\}_{z \in U_{t-1}}\) are independent. The second equality in (\ref{var_y}) holds since the random variable~\(X_{y,z}\) takes the value either~\(0\) or~\(1.\) The final estimate in (\ref{var_y}) follows from (\ref{mean_jy}).

Substituting the estimates~(\ref{mean_jy}) and (\ref{var_y}) into (\ref{tempax_123}) we get
\begin{eqnarray}
\mathbb{E}\left(\sum_{y \in S_{t-1}}J_y\right)^2 \leq C_u\#S_{t-1}   + (C_u\#S_{t-1})^2. \nonumber
\end{eqnarray}
and using the above in~(\ref{tempax_124}) we have
\begin{eqnarray}
\mathbb{E}(\#{\cal N}_t)^2\ind(F_{t-1}) \leq \left(C_u\#S_{t-1}   + C_u^2(\#S_{t-1})^2\right)\mathbb{P}(F_{t-1}). \nonumber
\end{eqnarray}
From the expression~(\ref{exp_nt2}), we then have
\begin{eqnarray}
\mathbb{E}(\#{\cal N}_t)^2 &\leq& \sum_{S_0 = \{1\},S_1,\ldots, S_{t-1}}\left(C_u\#S_{t-1}   + C_u^2(\#S_{t-1})^2\right)\mathbb{P}(F_{t-1})  \nonumber\\
&=& C_u \mathbb{E}(\#{\cal N}_{t-1}) + C_u^2\mathbb{E}(\#{\cal N}_{t-1})^2 \nonumber\\
&\leq& C_u^t + C_u^2\mathbb{E}(\#{\cal N}_{t-1})^2 \label{rec_eq111}
\end{eqnarray}
where the final estimate follows from (\ref{diam_exp}).

Setting \(a_t = \mathbb{E}(\#{\cal N}_t)^2\) and proceeding iteratively using (\ref{rec_eq111}), we have
\[a_t \leq \sum_{j=t}^{2t-1} C_u^{j} + C_u^{2t}a_0 = \sum_{j=t}^{2t} C_u^{j}\] since \(a_0 = \mathbb{E}\#({\cal N}_0)^2 = 1.\) Therefore we have that
\[a_t \leq \sum_{j\geq t} C_u^{j} = \frac{C_u^{t}}{1-C_u}.\]  This proves (\ref{diam_var}). \(\qed\)

\subsection*{Large components estimates}
Using Lemma~\ref{lem_diam}, we see that large components cannot exist with high probability.
\begin{Lemma}\label{larg_comp_sub} Fix \(C < 1.\) For every \(\epsilon > 0,\) there are positive constants \(D = D(C,\epsilon)\) and \(N = N(C,\epsilon)\) such that
\begin{eqnarray}
\mathbb{P}\left(\#{\cal E}_i \geq \epsilon n\right) \leq D \frac{(\log{n})^3}{n^2} \label{tau_use_imp}
\end{eqnarray}
for all \(n \geq N\) and  all \(1 \leq i \leq n.\)\\
\end{Lemma}
\emph{Proof of Lemma~\ref{larg_comp_sub}}: As in the proof of Lemma~\ref{lem_diam}, we set \(i = 1\) throughout and let \({\cal N}_t\)  denote the (random) set of vertices at distance \(t \geq 1\) from vertex \(1.\) The proof holds for any \(1 \leq i \leq t.\) Let
\begin{equation}\label{tau_def}
\tau = \inf\{{t \geq 1} : {\cal N}_t = \emptyset\}
\end{equation}
denote the largest distance of a vertex from the vertex~\(1.\) Fixing \(C < C_1 < 1,\) we have that \(C_u = C + \alpha_n \leq C_1\) for all \(n \geq n_0(C,C_1).\) The above statement is true since \(\alpha_n \longrightarrow 0\) as \(n \rightarrow \infty.\) Using (\ref{diam_exp}), we have for \(t = \frac{-2}{\log{C_1}} \log{n}\)  that
\begin{equation}\label{tau_max}
\mathbb{P}\left(\tau \geq t\right) \leq \mathbb{P}\left(\#{\cal N}_t \geq 1\right) \leq \mathbb{E}\#{\cal N}_t \leq C_u^{t}= \exp\left(-2\frac{\log{C_u}}{\log{C_1}}\log{n}\right) \leq \frac{1}{n^2}
\end{equation}
where the final inequality follows since \(C_u  \leq C_1 <1\) and so \(\frac{\log{C_u}}{\log{C_1}} \geq 1.\)
Thus for \(\epsilon > 0,\) we have
\begin{eqnarray}
\mathbb{P}\left(\#{\cal E}_1 \geq \epsilon n\right) &=& \mathbb{P}\left(\{\#{\cal E}_1 \geq \epsilon n\} \cap \{\tau \geq t\} \right) + \mathbb{P}\left(\{\#{\cal E}_1 \geq \epsilon n\} \cap \{\tau < t\} \right)\nonumber\\
&\leq& \mathbb{P}\left(\{\tau \geq t\} \right) + \mathbb{P}\left(\{\#{\cal E}_1 \geq \epsilon n\} \cap \{\tau < t\} \right)\nonumber\\
&\leq& \frac{1}{n^2} + \mathbb{P}\left(\{\#{\cal E}_1 \geq \epsilon n\} \cap \{\tau < t\} \right) \label{tau_use1}
\end{eqnarray}
using (\ref{tau_max}).

To evaluate the second term in (\ref{tau_use1}), we write
\begin{eqnarray}\label{tau_use2}
\mathbb{P}\left(\{\#{\cal E}_1 \geq \epsilon n\}\cap \{\tau < t\}\right) &=& \sum_{k=1}^{t-1} \mathbb{P}\left(\{\#{\cal E}_1 \geq \epsilon n\}\cap \{\tau = k\}\right)
\end{eqnarray}
and fix \(1 \leq k \leq t-1.\)
If the event~\(\{\#{\cal E}_1 \geq \epsilon n\}\cap \{\tau = k\}\) occurs, then some \({\cal N}_j, 1 \leq j \leq k\) contains at least \(\frac{\epsilon n}{k} \geq \frac{\epsilon n}{t}\) vertices and so
\begin{equation}
\mathbb{P}\left(\{\#{\cal E}_1 \geq \epsilon n\}\cap \{\tau = k\}\right) \leq \mathbb{P}\left(\bigcup_{j=1}^{k} \{\#{\cal N}_j \geq \epsilon n t^{-1}\} \right)
\leq \sum_{j=1}^{k} \mathbb{P}\left(\#{\cal N}_j \geq \epsilon n t^{-1}\right).\label{tau_use3}
\end{equation}
For any fixed \(1 \leq j \leq k,\) we have that
\begin{equation}\label{tau_use4}
\mathbb{P}\left(\#{\cal N}_j \geq \epsilon n t^{-1}\right)  \leq \frac{t^2}{(\epsilon n)^2}\mathbb{E}(\#{\cal N}_j)^2 \leq \frac{t^2}{(\epsilon n)^2}\frac{C_u^j}{1-C_u}
\end{equation}
where the first inequality follows using Markov inequality and the second inequality follows from the estimate~(\ref{diam_var}).

Using (\ref{tau_use4}), the final term in (\ref{tau_use3}) can be bounded above by
\begin{equation}\label{tau_use5}
\sum_{j=1}^{k} \frac{t^2}{(\epsilon n)^2}\frac{C_u^j}{1-C_u} \leq \sum_{j\geq 1 } \frac{t^2}{(\epsilon n)^2}\frac{C_u^j}{1-C_u} = \frac{t^2}{\epsilon^2n^2} \frac{C_u}{(1-C_u)^2} \leq D_1 \frac{t^2}{n^2},
\end{equation}
for all \(n \geq n_1\) and some constant \(D_1 = D_1(C,\epsilon) > 0.\) Here \(n_1\) does not depend on the choice of \(k.\) The final estimate holds since \(C_u = C+\alpha_n \longrightarrow C < 1\) as \(n \rightarrow \infty\) and so \(\frac{C_u}{(1-C_u)^2}\) is a bounded sequence in \(n.\)
The final term in (\ref{tau_use5}) is an estimate for \(\mathbb{P}\left(\{\#{\cal E}_1 \geq \epsilon n\}\cap \{\tau = k\}\right) \) and holds for any \(1 \leq k \leq t-1.\)
Substituting the above estimate into (\ref{tau_use2}) we have
\begin{eqnarray}\label{tau_use33}
\mathbb{P}\left(\{\#{\cal E}_1 \geq \epsilon n\}\cap \{\tau < t\}\right) \leq tD_1\frac{t^2}{n^2} \leq D_2 \frac{(\log{n})^3}{n^2}
\end{eqnarray}
for some constant \(D_2 = D_2(C,C_1,\epsilon) > 0\) and all \(n \geq n_2\) large. The final estimate holds since \(t = -\frac{2}{\log{C_1}}\log{n}.\) Finally, substituting (\ref{tau_use33}) into (\ref{tau_use1}) we get (\ref{tau_use_imp}).\(\qed\)

\emph{Proof of Theorem~\ref{th2}}:
Fix \(\epsilon > 0\) small so that~\(\delta_0(C,\epsilon,\epsilon)\) defined in (\ref{delta0_def}) is positive and let~\(\delta_0= \delta_0(C,\epsilon,\epsilon).\)

Let
\begin{equation}\label{no_giant}
W(\epsilon) := \bigcup_{1 \leq i \leq n} \{\#{\cal E}_i \geq \epsilon n\}
\end{equation}
denote the event that there is a component of size larger than \(\epsilon n.\)
From (\ref{tau_use_imp}) of Lemma~\ref{larg_comp_sub}, we have that
\begin{equation}\label{w_est}
\mathbb{P}(W(\epsilon)) \leq \sum_{i=1}^{n}\mathbb{P}\left(\#{\cal E}_i \geq \epsilon n\right) \leq D\frac{(\log{n})^3}{n}
\end{equation}
where the constant \(D\) is as in (\ref{tau_use_imp}). This estimates that
large components cannot exist with high probability.

To see existence of mid size components, we use Lemma~\ref{lem_mid}.
Let \(\delta > 0\) be as in (\ref{del_def2}) and let \(M > \frac{1}{\delta}\) be fixed and fix
\(\theta  < \min(M\delta-1,1).\) Since \(\delta_0 = \delta - C\epsilon -\epsilon \longrightarrow \delta\)
as \(\epsilon \rightarrow 0,\) we fix \(\epsilon > 0\) small enough so that
\(M > \frac{1}{\delta_0}\) and \(\theta  <  \min(M\delta_0-1,1).\)
Fix \(\theta_1> 0\) so that \(\theta < \theta_1 < \min(M\delta_0-1,)\)
and let \(B(M,\epsilon)\) be the event defined in~(\ref{b_m_def}) that
there is a component whose size lies in the range \([M\log{n}+1,\epsilon n].\) From~(\ref{bm_est}) we have
\begin{equation}\label{bm_est44}
\mathbb{P}(B(M,\epsilon)) \leq \frac{1}{n^{\theta_1}}
\end{equation}
for all \(n\) large. 

If \(B^c(M,\epsilon) \cap W^c(\epsilon)\) occurs, then every component has size at most \(M \log{n}.\) In other words, the event \(H_1(M)\)
defined in (\ref{h_2_def}) occurs. From~(\ref{bm_est44}) and~(\ref{w_est}), we have that
\[\mathbb{P}(H_1) \geq 1 - \frac{1}{n^{\theta_1}}  - D_2\frac{(\log{n})^3}{n} \geq 1 - \frac{1}{n^{\theta}}\] for all
\(n\) large. This proves (\ref{th2}).\(\qed\)



\setcounter{equation}{0}
\renewcommand\theequation{\thesection.\arabic{equation}}
\section{Proof of Theorem~\ref{th_q0}}\label{pf_q0}
To prove Theorem~\ref{th_q0}, we need the following properties of the function \(q(C)\) defined in (\ref{q_0_def}).\\
\((a1)\) The function \(q(C)\) is continuous at any \(C> 0.\)\\
\((a2)\) If \(C_1 > C_2 \geq 1,\) then \(q(C_1) < q(C_2)\) strictly.\\\\
\emph{Proof of \((a1)-(a2)\)}: To prove the continuity property \((a1),\) we use an approximation procedure and truncate the series expansion for \(q(C)\) at a finite number of terms.
For \(C > 0,\) we let \[q_{0,N}(C) = \sum_{r=1}^{N} \frac{T_rC^{r-1}e^{-Cr}}{(r-1)!}\] and have
\(q_{0,N}(C) \uparrow q(C)\) as \(N \rightarrow \infty.\)
The function \(q_{0,N}(C)\) is continuous function of \(C\) for any fixed \(N.\)
Let \(C_0 > 0\) be fixed. We have
\begin{eqnarray}
|q(C_0) - q(C)| &\leq& |q(C_0) - q_{0,N}(C_0)| + |q_{0,N}(C_0) - q_{0,N}(C)| \nonumber\\
&&\;\;\;\;+ |q_{0,N}(C) - q(C)| \nonumber\\
&=& (q(C_0) - q_{0,N}(C_0)) + |q_{0,N}(C_0) - q_{0,N}(C)| \nonumber\\
&&\;\;\;\;+ (q_{0,N}(C) - q(C)) \nonumber\\
&=& R_N(C_0) + |q_{0,N}(C_0) - q_{0,N}(C)| + R_N(C) \nonumber
\end{eqnarray}
where the second equality follows from the fact that \(q_{0,N}(.) \leq q(.)\) by definition. The term \(R_N(.)\) is the remainder
term defined in (\ref{rem_def}).
We also note that for any \(C > 0,\) we have
\[R_N(C) = \sum_{r = N+1}\frac{T_re^{-r}}{(r-1)!}e^{-\delta r} \leq \sum_{r = N+1}\frac{T_re^{-r}}{(r-1)!} = R_N(1).\]
where \(\delta = \delta(C) \geq 0\) is as in (\ref{del_def2}). 
Thus
\begin{eqnarray}\label{q0_trunc_est}
|q(C_0) - q(C)| \leq |q_{0,N}(C_0) - q_{0,N}(C)| + 2R_N(1)
\end{eqnarray}

Fix \(C_0 > 0\) and let \(\{C_m\}_{m \geq 1}\) be any sequence such that \(C_m > 0\) for all \(m \geq 1\) and \(C_m \longrightarrow C_0\) as \(m \rightarrow \infty.\) Fix integer \(N \geq 1.\) Using (\ref{q0_trunc_est})
and the continuity of \(q_{0,N}\) we have
\[\limsup_m |q(C_0) - q(C_m)| \leq 2R_N(1).\] Allowing \(N \rightarrow \infty\) and using (\ref{rem_est}) then gives \(\lim_m q(C_m) = q(C_0).\) This proves that \(q(C)\) is a continuous function for \(C > 0.\)

To see the strictly decreasing property \((a2),\) we fix \(C_1 > C_2 \geq 1.\) We then have \(C_1e^{-C_1} < C_2 e^{-C_2}\) and \(\frac{1}{C_1} < \frac{1}{C_2}.\) Therefore for any integer \(r \geq 1,\) we have \[C_1^{r-1}e^{-C_1r} = \frac{1}{C_1} (C_1e^{-C_1})^r  < \frac{1}{C_2} (C_2e^{-C_2})^r = C_2^{r-1}e^{-C_2 r}.\] Multiplying by \(\frac{T_r}{(r-1)!}\) both sides and summing over \(r\) gives that \(q(C_1) < q(C_2)\) strictly. This proves the strictly decreasing property of~\(q(C)\) for \(C > 1.\) \(\qed\)\\\\
\emph{Proof of Theorem~\ref{th_q0}}: \((i)\) Proved before, see~(\ref{qc_less_one}). \\
\((ii)\) This is true by the property of
the extinction probability for the Poisson branching process (see Durrett (2007)).
For completeness, however, we give a small proof directly using the properties of random graphs.

Fix \(0 < \gamma < 1\) and \(C \neq 1.\) Let \(\epsilon > 0\) be such that (\ref{exp_x_est}) holds
and let \(Z_n(\epsilon)\) be as defined in (\ref{x_def}). We have from~(\ref{exp_x_est})
that \[\mathbb{E}Z_n(\epsilon) \leq nq(C)(1+\gamma)\] for all \(n\) large. We also have that \(Z_n(\epsilon) = n\) if and only if there is no component of size
larger than \(\epsilon n\) i.e., if and only if the event \(W^c(\epsilon)\) holds.
Here~\(W(\epsilon)\) is as defined in (\ref{no_giant}). Using the corresponding estimate (\ref{w_est}) we have
that \[nq(C)(1+ \gamma) \geq \mathbb{E}Z_n(\epsilon) \geq n\mathbb{P}(W^c(\epsilon)) \geq n\left(1 - D_2\frac{(\log{n})^3}{n}\right)\]
where the constant \(D_2 > 0\) is as in (\ref{w_est}). Therefore \[q(C) \geq \frac{1}{1+\gamma}\left(1 - D_2\frac{(\log{n})^3}{n}\right)\] for all \(n\) large.
Allowing \(n \rightarrow \infty\) and then \(\gamma \downarrow 0\) gives \(q(C) \geq 1.\)
We have so far proved that \(q(C) = 1\) for \(0 <C < 1.\)

To prove \((iii)\) and that \(q(1) = 1,\) we use the properties \((a1)-(a2)\) of \(q(.)\) described above.

\((iii)\) Since \(q(C)\) is continuous at \(C = 1\) (see property \((a1)\)), we have
\[\lim_{C \uparrow 1} q(C) = q(1) = 1.\]
Also we have that \(q(C)\) is strictly decreasing in \(C\) for \(C > 1\) (see property \((a2)\))
and so
\begin{equation}\label{q_ineq2}
0 < q(C) < q(1) = 1
\end{equation}
for all \(C > 1.\) This proves \((iii).\)

\((iv)\) The proof is combinatorial and we refer to Chapter~\(14,\) van Lint and Wilson~(2010).\(\qed\)

\setcounter{equation}{0}
\renewcommand\theequation{\thesection.\arabic{equation}}
\section{Supercritical case}\label{pf_sup}





\emph{Proof of Theorem~\ref{th1}}: We recall that \(Z_n(\epsilon)\) defined in (\ref{x_def}) denotes the sum of sizes of components each of whose size is at most \(\epsilon n.\) We use the estimates on mean and variance of~\(Z_n(\epsilon)\) in Lemma~\ref{exp_x_lem} and Lemma~\ref{var_x_lem} to see that at least one giant component exists.

Fix \(0 < \gamma < 1, C > 1\) and let \(\delta > 0\) be as in (\ref{del_def2}). Define \(\epsilon_3  = \min\left(\epsilon_1,\epsilon_2\right)\) where \(\epsilon_1\) and \(\epsilon_2\) are defined in Lemma~\ref{exp_x_lem} and~\ref{var_x_lem}, respectively. Let \(q = q(C)\) be as in~(\ref{q_0_def}).

Using the Chebychev's inequality we have
\begin{equation}
\mathbb{P}\left(|Z_n(\epsilon)-\mathbb{E}Z_n(\epsilon)| > q\frac{\gamma}{2} n\right) \leq \frac{4var(Z_n(\epsilon))}{q^2\gamma^2n^2} \label{var_zn_cheb}
\end{equation}
where \(var(Z_n(\epsilon)) = \mathbb{E}(Z_n(\epsilon)-\mathbb{E}Z_n(\epsilon))^2\) is the variance of the random variable~\(Z_n(\epsilon).\)
Since \(0 < \gamma < 1\) is arbitrary, we use Lemma~\ref{var_x_lem} with \(\gamma_1 = \frac{\gamma^3}{64}\) to get that \(var(Z_n(\epsilon)) \leq nq\left(1+\frac{\gamma^3}{64}\right) + n^2q^2 \frac{\gamma^3}{16}\) for all \(n\) large. Therefore the final term in (\ref{var_zn_cheb}) is bounded above by
\begin{equation}
\frac{4}{q^2\gamma^{2}n^2}\left(nq\left(1+\frac{\gamma^3}{64}\right) + n^2q^2\frac{\gamma^3}{16}\right) = \frac{4}{q\gamma^{2}n}\left(1+\frac{\gamma^3}{64}\right) +  \frac{\gamma}{4} \leq \frac{\gamma}{2} \label{v_eps_def}
\end{equation}
for all \(n\) large.

Thus if \(Y_n(\epsilon) = n-Z_n(\epsilon)\) denotes the sum of sizes of components each of which has size at least \(\epsilon n\) as defined in (\ref{y_def}), we have from (\ref{var_zn_cheb}) and~(\ref{v_eps_def}) that
\begin{eqnarray}
\mathbb{P}\left(|Y_n(\epsilon)-\mathbb{E}Y_n(\epsilon)| > q \frac{\gamma}{2}n\right) \leq \frac{\gamma}{2}\label{giant_dev_est}
\end{eqnarray}
for all \(n\) large. Since \(q = q(C) < q(1) = 1\) strictly for \(C > 1\) (see properties \((ii)-(iii)\) of Theorem~\ref{th_q0}), we have that there exists at least one giant component.
Using Lemma~\ref{exp_x_lem} with \(\frac{\gamma}{2},\) we also have
\begin{equation}\label{y_n_est11}
n\left(1-q-\frac{\gamma}{2}\right)  \leq \mathbb{E}Y_n(\epsilon) \leq n\left(1-q + \frac{\gamma}{2}\right)
\end{equation}
for all \(n\) large. If the estimate~(\ref{y_n_est11}) holds and the event on the left hand size of~(\ref{giant_dev_est}) holds, then the event~\(V(\gamma,\epsilon)\) defined in~(\ref{v_i_def}) occurs. We therefore have from~(\ref{giant_dev_est}) and~(\ref{y_n_est11}) that
\begin{equation}
\mathbb{P}(V(\gamma,\epsilon)) \geq 1 - \frac{\gamma}{2}\label{v_1_est}
\end{equation}
for all \(n\) large.

Next, we recall the event \(B(M,\epsilon)\) defined in (\ref{b_m_def}) that there exists a component in the range \([M\log{n}+1, \epsilon n].\) Let \(\delta > 0\) be as in (\ref{del_def2}). Fix \( M >\frac{1}{\delta}\) and choose \(\epsilon_3\) (defined in the first paragraph of this proof) smaller if necessary so that the following holds for all \(0 < \epsilon < \epsilon_3.\) If \(\delta_0 = \delta_0(C,\epsilon,\epsilon) = \delta - (2C+1)\epsilon,\) then \(\delta_0 > 0\) and \(M > \frac{1}{\delta_0}.\) Fixing \(0 < \epsilon < \epsilon_3\) and \(0 < \theta < \min(M\delta_0-1,1),\) we have from~(\ref{bm_est}) that
\begin{equation}\label{mid_comp111}
\mathbb{P}(B(M,\epsilon)) \leq \frac{1}{n^\theta}
\end{equation}
for all \(n\) large.

From (\ref{v_1_est}) and (\ref{mid_comp111}), we note that
\begin{equation}
\mathbb{P}\left(V^c(\gamma,\epsilon)\cup B(M,\epsilon)\right) \leq \frac{\gamma}{2} + \frac{1}{n^{\theta}} \leq \gamma \label{h_1_est111}
\end{equation}
for all \(n\) large.

If \(V(\gamma,\epsilon) \cap B^c(M,\epsilon)\) occurs, then there is at least one component containing at least \(\epsilon n\) nodes. Also, since \(B^c(M,\epsilon)\) occurs,  size of all small components is at most \(M\log{n}.\) And by definition, the event \(H_2 = H_2(M,\gamma,\epsilon)\) defined in~(\ref{h_1_def}) occurs. From (\ref{h_1_est111}), we obtain Theorem~\ref{th1}.\(\qed\)

\setcounter{equation}{0}
\renewcommand\theequation{\thesection.\arabic{equation}}
\section{Proof of Theorems~\ref{th111} and~\ref{th112}} \label{pf11}
\emph{Proof of Theorem~\ref{th111}}: We first see that there is a unique solution to \(\delta_0\left(C,\frac{1}{2},0\right) = \frac{C}{2} - 1 - \log{C} = 0\) for \(C > 1.\) The derivative \(\delta'_0 = \delta'_0(C,\frac{1}{2},0)\) with respect to \(C,\) satisfies \(\delta'_0 = \frac{1}{2} - \frac{1}{C} > 0\) for \(C > 2\) and so the function \(\delta_0\left(C,\frac{1}{2},0\right)\) is increasing for all \(C > 2.\) Moreover, \(\delta_0\left(2,\frac{1}{2},0\right) < 0\) and \(\delta_0\left(C,\frac{1}{2},0\right) \longrightarrow \infty\) as \(C \rightarrow \infty.\) Therefore there is a unique \(C_0 > 2\) such that \(\delta_0\left(C_0,\frac{1}{2},0\right) = 0.\)

To see that there is no solution for \(\delta_0\left(C,\frac{1}{2},0\right) = 0\) in the interval \((1,2),\) we use the fact that the derivative \(\delta'_0  = \frac{1}{2} - \frac{1}{C} < 0\) for all \(C < 2.\) Thus \(\delta_0\left(C,\frac{1}{2},0\right)\) is decreasing in the interval \((1,2)\) and \(\delta_0\left(1,\frac{1}{2},0\right) = -\frac{1}{2} < 0.\) This implies that there is no solution for \(\delta_0\left(C,\frac{1}{2},0\right) = 0\) in the interval \((1,2).\) We also have \(C_0 < 8,\) since \(\delta_0\left(8,\frac{1}{2},0\right) = 3 - \log{8}  = 3(1-\log{2}) > 0.\)

To obtain the uniqueness of the giant component, we use the upper bound in the tree counting estimate (\ref{c1_est1}) of Lemma~\ref{c1_est1_lem}. To ensure that the upper bound is valid, we fix \(C > C_0\) and let \(\epsilon_0 > \frac{1}{2}\) and \(\alpha,\omega_0 > 0\) be such that
\begin{equation}\label{del_choice}
0 < \delta_0\left(C,\alpha+\epsilon_0,\omega_0\right) < \delta_0\left(C,\epsilon_0,\omega_0\right)
\end{equation}
If \(Z_n(\epsilon_0)\) as defined in (\ref{x_def}) is the number of vertices belonging to non giant components i.e., components with size at most~\(\epsilon_0 n,\)
we have from (\ref{c1_est1}) that
\begin{equation}
\frac{1}{n}\mathbb{E}Z_n(\epsilon_0) = \frac{1}{n}\sum_{i=1}^{n}\mathbb{P}(1 \leq \#{\cal E}_i \leq \epsilon_0 n) \leq \sum_{r=1}^{\epsilon_0 n} \frac{1}{C}\frac{T_re^{-r}}{(r-1)!}e^{-\delta_0 r} \leq q_0(C,\epsilon_0,\omega_0).\label{exp_x_est22}
\end{equation}
Here \(\delta_0 = \delta_0(C,\epsilon_0,\omega_0)\) and \(q_0 := q_0(C,\epsilon_0,\omega_0)\) is as defined in (\ref{q_i_def}).
Using Markov inequality, we have that
\begin{equation}\label{eq_sum}
\mathbb{P}(Z_n(\epsilon_0) > n e^{\alpha C} q_0) \leq \frac{\mathbb{E}Z_n(\epsilon_0)}{e^{\alpha C}q_0} \leq e^{-\alpha C}
\end{equation}
where \(\alpha > 0\) is as in (\ref{del_choice}).

We have that \(e^{\alpha C} q_0 < 1 \) strictly. To see this is true, we first have from the choice of \(\alpha > 0\) in (\ref{del_choice})
that \(\delta_{\alpha} := \delta_0(C,\epsilon_0+\alpha,\omega_0) > 0.\) Also from definition of \(q(C)\) in (\ref{q_0_def}) and from Theorem~\ref{th_q0} we have that \(q(1) = \sum_{r=1}^{\infty} \frac{T_re^{-r}}{(r-1)!} = 1.\) Therefore
\[q_0(C,\epsilon_0+\alpha,\omega_0) = \sum_{r=1}^{\infty} \frac{1}{C}\frac{T_re^{-r}}{(r-1)!}e^{-\delta_{\alpha} r}  \leq \frac{1}{C}\sum_{r=1}^{\infty} \frac{T_re^{-r}}{(r-1)!} = \frac{1}{C}  <1\] since \(C > 1.\) Since \(e^{-\delta_{\alpha} r} = e^{-\delta_0 r}e^{C\alpha r} \geq e^{-\delta_0 r}e^{C\alpha}\) for integer \(r \geq 1,\) we also have
\[q_0(C,\epsilon_0+\alpha,\omega) = \sum_{r=1}^{\infty} \frac{1}{C}\frac{T_re^{-r}}{(r-1)!}e^{-\delta_\alpha r}\geq e^{\alpha C}\sum_{r=1}^{\infty} \frac{1}{C}\frac{T_re^{-r}}{(r-1)!}e^{-\delta_0 r} = e^{\alpha C}q_0.\] Thus \(e^{\alpha C}q_0 < 1,\) strictly.

The term \(Z_n(\epsilon_0)\) is the sum of sizes of non giant components i.e., the components whose size lies between \(1\) and~\(\epsilon_0 n.\) If \(\{Z_n(\epsilon_0) \leq ne^{\alpha C}q_0\}\) occurs, then \(Z_n(\epsilon_0) < n\) strictly and therefore there exists a component with at least \(\epsilon_0 n +1 > \frac{n}{2} + 1\) vertices. This proves the first step in the proof.


The next step is to eliminate mid size components. Since the term \(\delta_0 = \delta_0(C,\epsilon_0,\omega_0) > 0,\) the result of Lemma~\ref{lem_mid} holds. In particular, for a fixed integer \(M_0 > \frac{1}{\delta_0}\) we have from estimate~(\ref{bm_est}) of Lemma~\ref{lem_mid} that
\begin{equation}\label{eq_b}
\mathbb{P}(B(M_0,\epsilon_0)) \leq \frac{1}{n^{\theta}}
\end{equation}
for some constant \(\theta > 0.\) Here \(B(M_0,\epsilon_0)\) denotes the event that there exists a component whose size lies in the range
\([M_0\log{n}+1,\epsilon_0 n].\)

From (\ref{eq_sum}) and (\ref{eq_b}) , we obtain that
\begin{equation}\label{ab_est}
\mathbb{P}\left(\{Z_n(\epsilon_0) \leq ne^{\alpha C}q_0\} \cap B^c(M_0,\epsilon_0)\right) \geq 1 - e^{-\alpha C}  - \frac{1}{n^{\theta}}
\end{equation}
for all \(n \geq n_0.\) Here \(n_0 \geq 1\) is a constant. If the event on the left hand side occurs, there is at least one component containing at least \(\epsilon_0 n > \frac{n}{2} + 10\) vertices, for all \(n\) large. This means that every other component has at most \(\frac{n}{2}-10\) vertices. But since \(B^c(M_0,\epsilon_0)\) also occurs, every other component has size at most \(M_0\log{n}.\) Therefore there is exactly one (giant) component containing at least \(\frac{n}{2} + 10\) vertices and every other component has at most~\(M_0\log{n}\) vertices.


The final step is to use the results of Theorem~\ref{th1} to obtain a better estimate for the size of the unique giant component. For \(0 < \gamma < 1,\epsilon > 0,\) let \(H_2 = H_2(M,\gamma,\epsilon)\) be as defined in the statement of Theorem~\ref{th1}. Fixing \(M > \frac{1}{\delta}\) and choosing \(\epsilon > 0\) small we have from (\ref{h_2_est}) of Theorem~\ref{th1} that \(\mathbb{P}(H_2) \geq 1 - \frac{\gamma}{2}.\) Using estimate (\ref{ab_est}), we therefore have that
\begin{equation}\label{ab_est2}
\mathbb{P}\left(\{Z_n(\epsilon_0) \leq ne^{\alpha C}q_0\} \cap B^c(M_0,\epsilon_0) \cap H_2 \right) \geq 1 - e^{-\alpha C}  - \frac{1}{n^{\theta}} - \frac{\gamma}{2} \geq 1-e^{-\alpha C} -\gamma
\end{equation}
for all \(n\) large.

Suppose that the event on the left hand side of (\ref{ab_est2}) occurs. By the discussion following (\ref{ab_est}), there is exactly one giant component containing at least \(\frac{n}{2} + 1\) nodes and the rest of the components have size at most \(M_0\log{n}.\) But since \(H_2\) also occurs, the size of the giant component is in the range \([(1-q(C)-\gamma)n,(1-q(C)+\gamma)n]\) and every other component has size at most \(M\log{n}.\) Here \(q(C)\) is as defined in (\ref{q_0_def}). This proves Theorem~\ref{th111}.~\(\qed\)

\emph{Proof of Theorem~\ref{th112}}:
If \(C < e^{-1},\) then the term \[\delta_0\left(C,\epsilon,0\right) = C(1-\epsilon) - 1 - \log{C} = C(1-\epsilon) + \log\left(\frac{e^{-1}}{C}\right) > 0\] for all \(0 < \epsilon \leq 1.\) Fixing \(\epsilon = 1,\) we have that \(\delta_0 = \delta_0(C,1,0) = \log\left(\frac{e^{-1}}{C}\right) > 0.\) Therefore the upper bound (\ref{c1_est1}) of Lemma~\ref{c1_est1_lem} holds for all \(1 \leq r \leq n.\) Fixing \( M > 0,\) we therefore have that
\begin{equation}\label{larg_est}
\mathbb{P}(\#{\cal E}_i \geq M\log{n}) \leq \sum_{r = M\log{n}}^{n} \frac{T_re^{-r}}{C(r-1)!} e^{-\delta_0r} \leq e^{-\delta_0M\log{n}}\sum_{r=M\log{n}}^{n}\frac{T_re^{-r}}{C(r-1)!} .
\end{equation}
The summation in the final term is bounded above as \[\sum_{r=M\log{n}}^{n}\frac{T_re^{-r}}{C(r-1)!} \leq \frac{1}{C} \sum_{r=1}^{\infty} \frac{T_re^{-r}}{(r-1)!} = \frac{q(1)}{C} = \frac{1}{C}\] using Theorem~\ref{th_q0}. Substituting in (\ref{larg_est}) we have
\begin{equation}\nonumber
\mathbb{P}(\#{\cal E}_i \geq M\log{n}) \leq \frac{1}{C} e^{-\delta_0M\log{n}} = \frac{1}{C} \frac{1}{n^{M\delta_0}} \leq \frac{1}{n^{\theta+1}}
\end{equation}
for all \(n\) large, provided \(M > 0\) is large. Thus
\begin{equation}\nonumber
\mathbb{P}\left(\bigcup_{1 \leq i \leq n}\#{\cal E}_i \geq M\log{n}\right) \leq \frac{1}{n^{\theta}}.
\end{equation}
This implies that with probability at least \(1-\frac{1}{n^{\theta}}\) all the components have size at most \(M\log{n}.\)\(\qed\)

\section*{Acknowledgements}
I thank Professors Rahul Roy, Thomas Mountford, Siva Athreya and Federico Camia for crucial comments and for my fellowships. I also thank ISF-UGC for my fellowship.



\setcounter{equation}{0} \setcounter{Lemma}{0} \renewcommand{\theLemma}{II.%
\arabic{Lemma}} \renewcommand{\theequation}{II.\arabic{equation}} %
\setlength{\parindent}{0pt}




%





\bibliographystyle{plain}

\end{document}